\documentclass[12pt,a4paper,reqno]{amsart}

\usepackage{amsmath,amsthm,amssymb}
\usepackage[english]{babel}
\usepackage[dvips]{graphicx}
\usepackage{hyperref}

\numberwithin{equation}{section}

\newtheorem{proposition}{Proposition}

\newtheorem{definition}{Definition}

\newcommand {\bs} {\boldsymbol}
\newcommand {\A} {\mathop{\rm Acc}\nolimits}
\newcommand {\rk} {\mathop{\rm rank}\nolimits}
\newcommand {\ri} {\mathrm{i}}
\newcommand {\diag} {\mathop{\rm diag}\nolimits}
\newcommand {\mbf}[1] {\mathbf{#1}}
\newcommand {\mF}{\mathcal{F}}
\newcommand {\mP}{\mathcal{P}}
\newcommand {\bR}{\mathbf{R}}
\newcommand {\mm}{\mathcal{M}_1}
\newcommand {\mn}{\mathcal{M}_2}
\newcommand {\mo}{\mathcal{M}_3}
\newcommand {\ds}{\displaystyle}
\newcommand {\mstrut}{\vphantom{\bigl(}}

\begin{document}

\title[Geometrical approach to separation of variables]{Geometrical approach to separation of variables\\in mechanical systems}

\author{Mikhail P. Kharlamov, Alexander Y. Savushkin}
\address{Russian Presidential Academy of National Economy and Public Administration, Volgograd branch}
\email{mharlamov@vags.ru}

\date{20.05.2010}

\begin{abstract}
The article presents a compact review of the analytical results (2002--2009) in the study of the system describing the motion of a top in two constant fields. The Liouville integrability of this system under certain condition of the Kowalevski type was established by A.G.\,Reyman and M.A.\,Semenov-Tian-Shansky. We present some geometrical foundations of finding separations of variables. Two systems of local planar coordinates are introduced leading to separation of variables for two subsystems with two degrees of freedom in the dynamics of the generalized Kowalevski top.
\end{abstract}


\maketitle

\tableofcontents

\section{Introduction}\label{sec1}

Historically, the main goal in the study of a differential equation of mechanics was to integrate it, to find the explicit solution. The fact of the existence of principally non-integrable systems has shifted the accent to qualitative methods of investigation. However, during the last 20 years new algebraic approaches to construct or reveal the integrable systems were created. For Hamiltonian systems, many new cases of commutative and non-commutative integrability were found. A lot of them deal with the problem of the rotation of a rigid body about a fixed point and its generalizations to Lie groups that differ from $SO(3)$. Here lies one of the reasons for a new wave of mathematical interest in the integrable Hamiltonian systems. Note the books \cite{BorMam0,BorMam,BorMam2,ReySemBk,TsiBk} published lately and containing more than a thousand of references on the matter. In relatively recent studies, the integrable systems with three or more degrees of freedom were found not having continuous symmetries and therefore not reducible as a whole to families of systems with two degrees of freedom (the so-called irreducible systems). One of the most ingenious and complicated example of the integrable irreducible system is the generalized Kowalevski gyrostat. It is a rigid body with a rotor satisfying the conditions of the Kowalevski type and placed in two independent constant fields, say, gravitational and magnetic. This case was found due to the efforts of several mathematicians. We note the works of O.I.\,Bogoyavlensky \cite{BogRus1,BogRus2,BogEn}, H.M.\,Yehia \cite{Yeh1,Yeh2}, L.N.\,Gavrilov \cite{Gavr}, I.V.\,Komarov \cite{Komar}. The final and mostly general result was obtained by A.G.\,Reyman and M.A.\,Semenov-Tian-Shansky \cite{ReySem,ReySemRus}. Nevertheless, it appeared that no explicit solutions in the classical sense can be found for the general case. Up-to-date the only progress in studying the Reyman--Semenov-Tian-Shansky case is based on the idea of investigating critical subsystems. This row of publications starts from the work \cite{ZotRCD}, in which the Bogoyavlensky system generalizing the 1st Appelrot class is investigated without any connection to the general case, and the work \cite{Odin}, in which the next generalization, of the 2nd and 3d Appelrot classes, was found. Later, the bifurcation diagrams of the partial cases and of the general case including the case of a gyrostat were built and classified in the space of parameters and with respect to the energy levels in the series of works \cite{KhSavSh33,Kh34,KhSh34,KhRCD05,Kh35,KhShRCD06,KhND07}. The complete description was obtained for the critical set of the momentum mapping stratified according to various criteria \cite{Kh362}. For two critical subsystems generalizing the famous Appelrot classes separations of variables were found and the phase topology was investigated in \cite{KhSavDan,KhSavMRC,KhSav, KhND06,Kh38,KhRCD09}, all critical periodic motions were described analytically and integrated in elliptic functions for the cases of a top \cite{Kh361,Kh34} and of a gyrostat \cite{Kh37,KhHMJ}. From the integrability point of view, each of the critical subsystems found is as general as the whole Kowalevski case in the gravity field, i.e., generates a one-parametric family of Hamiltonian systems with two degrees of freedom\footnote{In \cite{KhHMJ}, all critical subsystems for the Kowalevski gyrostat in two constant fields are found having one extra parameter.}. These systems have all the features of the systems in which separation of variables is obtained algebraically. Here we call such systems algebraically solvable. We present though the {\it geometrical} ideas of how to construct such separations.

A lot of publications deal with the notion of an algebraic complete integrable systems. This notion does not have a unique and exact definition. Usually it is associated with Hamiltonian systems the topology of which is connected with the Jacobians and Prim varieties of algebraic surfaces. Algebraic complete integrable systems with known separation of variables usually also have algebraic expressions of the initial phase variables in terms of the separated ones, although many contemporary papers, not like the classical works, do not pay much attention to this point. We think that a system with separated variables cannot be considered completely analytically integrated without explicit expressions for the phase variables. Separations of variables on algebraic curves generated by Lax representations, the existence of which is theoretically proved, do not provide any method to establish the connection with the initial variables. Here we use another approach --- starting from the geometry of initial variables and trying to find such projections of the phase space which bring the integral manifolds to the most simple image. If it is done (and there are no general ways to do it), then the connection with the initial variables is kept from the beginning. In this paper we propose the projections and the coordinate systems that made it possible to find new algebraically and topologically interesting separations of variables. The examples given below show the valuable applications.

\section{Geometrical approach to separation of variables}\label{sec2}
\subsection{Geometry of separated variables}
The term we use here---algebraically solvable systems---is not a standard one. Let us give some explanations of the problem. Consider a system of differential equations on a submanifold $\mathcal{P}$ in a real vector space
\begin{equation}
\label{neq2_1} \frac{d \mathbf{x}}{dt} = \mathbf{X}(\mathbf{x}),\quad \mathbf{x} \in \mathcal{P} \subset V
\end{equation}
defining the evolution of the set of phase variables ${\mbf{x}}$ and admitting separation of variables
\begin{equation}
\label{neq2_2} \frac {ds_i } {d\tau } = \sqrt {\mstrut P_i (s_i;\mbf{f})},
\end{equation}
where $s_i$ are so-called auxiliary or separated variables. We denote the vector of these variables by $\mathbf{s}$. The functions $P_i(s; \mathbf{f})$ are supposed to be polynomials in one variable $s$ with coefficients depending on a set of arbitrary constants ${\mbf{f}}$, $\tau $ is the ``reduced time'' connected with the real time $t$ by the equation of the type $d\tau / dt = \rho(\mbf{s})>0$. Usually such separation has a mechanical origin, therefore  equations (\ref{neq2_2}) have the form of an energy integral. Moreover, we suppose that all phase variables $x_j$ are expressed via the separated variables by rational functions with respect to the set of radicals
\begin{equation}
\notag R_{i\alpha }  = \sqrt {\mstrut s_i  - e_\alpha  } \quad \quad (e_\alpha \in {\mbf{C}})
\end{equation}
with coefficients smoothly depending on $\mathbf{s}$. The set of numbers $\{ e_\alpha\}$ depending, of course, on the constants ${\mbf{f}}$ includes also the roots of $P_i$. The classical examples show that the roots of $P_i$ usually provide the whole set of~$\{e_\alpha\}$.

The vector ${\mbf{f}}$ for the system (\ref{neq2_1}) usually appears as a collection of arbitrary constants of some set of first integrals $\{{\mathcal F}_j\}$. The dependencies
\begin{equation}
\label{neq2_4} {\mbf{x}} = {\mbf{x}}({\mbf{s}};{\mbf{f}}),\quad {\mbf{s}} \in \A({\mbf{f}})
\end{equation}
then are considered as the multi-valued parametric equations of an integral manifold $\mathcal{F}_\mbf{f} \doteq {\mathcal F}^{-1}({\mbf{f}}) \subset {\mathcal P}$. Here $\A({\mbf{f}})$ is the region in the separated variables space filled with trajectories of (\ref{neq2_2}) for given~${\mbf{f}}$. In Russian, it is traditionally called by the term corresponding to ``the region of the possibility of motions'' \cite{Tat1,Kh82}, which sounds not very elegant. A shorter term ``accessible region'' was introduced in English translations and completely reflects the situation. Notation of the region $\A({\mbf{f}})$ for the argument vector $\mbf{s}$ in (\ref{neq2_4}) is based on this term.

Let the system (\ref{neq2_1}) be a Hamiltonian complete integrable system with two degrees of freedom and
\begin{equation}
\label{neq2_5} \mathcal{F}:\mathcal{P} \to {\mbf{R}}^2
\end{equation}
be its integral mapping of compact character (the inverse image of a compact set is compact). Two degrees of freedom are chosen to simplify notation and because of the main examples below\footnote{In fact, here we do not use even the Hamiltonian type of the systems considered.}.

\begin{definition}\label{thdef21} The set $\operatorname{Im} \mathcal{F}$  is called the admissible region ${\rm (}$in the space of integral constants${\rm )}.$ A point ${\mbf{f}} \in {\bR}^2 $ is called an admissible point if and only if  ${\mathcal{F}^{ - 1} (\mbf{f}) \ne \varnothing}$. \end{definition}

Generally speaking, separation of variables is some mapping
\begin{equation}
\label{neq2_6} \pi :\mathcal{P}  \to {\mbf{R}}^2 _{(s_1 ,s_2 )}
\end{equation}
of the phase space to the plane of the auxiliary variables $(s_1 ,s_2 )$ that takes the system (\ref{neq2_1}) to the system (\ref{neq2_2}) with $i = 1,2$. For any admissible ${\mbf{f}}$ we obtain the mapping
\begin{equation}
\label{neq2_7} \pi _{\mbf{f}} :\mathcal{F}_{\mbf{f}}  \to {\mbf{R}}^2 _{(s_1 ,s_2 )}.
\end{equation}

\begin{definition}\label{thdef22} For a given ${\mbf{f}}$ the set $\A({\mbf{f}}) = \pi _{\mbf{f}} (\mathcal{F}_{\mbf{f}} )$ is called an accessible region in the plane of the separated variables. \end{definition} It is clear that $\A({\mbf{f}})$ is a subset of
\begin{equation}
\label{neq2_8} \{ (s_1 ,s_2 ):P_i (s_i ,{\mbf{f}}) \geqslant 0,i = 1,2\},
\end{equation}
and any point is included in $\A({\mbf{f}})$ with its total connected component in the set (\ref{neq2_8}). Then $\A({\mbf{f}})$ is a collection of rectangular regions in the $(s_1,s_2)$-plane;  ``infinite rectangles'' (half-strips, strips or even quadrants) are also possible. The fact that the accessible region is rectangular for all integral constants, obviously, is the necessary condition of separation of variables. Possibly, supposing that the projection of the type (\ref{neq2_7}) is non-trivial in some sense or that the boundaries of accessible regions essentially depend on $\mbf{f}$ in terms of the corresponding Jacobians, one can prove that this fact is also sufficient, but here we do not consider such a goal. We want to understand in which way to look for separations of variables.

Note that separation of variables not always has the form of a global mapping (\ref{neq2_6}). As far as the dynamics of a rigid body is concerned, separations not depending on integral constants are known for the Euler case (one independent variable can also be treated as a separation), the Goryachev--Chaplygin case and its generalization to a gyrostat given by L.N.\,Sretensky \cite{Sret}, the Clebsch case with zero momentum constant (isomorphic to the Neumann problem). Already in the Kowa\-lev\-ski case \cite{Kowa}, the change of variables leading to separation essentially uses the integral equations. However, having the mapping (\ref{neq2_7}) for \textit{all} $\mbf{f}$ one can formally substitute $\mbf{f}=\mF(\mbf{x})$ to obtain some mapping of the type (\ref{neq2_6}). It is important to mention that practically in all cases separation of variables is not built straightforwardly in the form (\ref{neq2_6}). As a rule, some intermediate mapping $\mP \to \bR^2_{\bs \xi}$ is constructed first. Let us call this auxiliary plane of the variables $\xi_1,\xi_2$ {\it the picture plane}. With the projection onto the picture plane, the separated variables are introduced as the functions of $\xi_1,\xi_2,\mbf{f}$ (usually defined by some implicit equations).

What is the geometry behind this construction? The mapping of the phase space to the picture plane takes an integral manifold to a region $\mF_{\mbf{f}}$ bounded by some curves. If these curves could be included in a two-parametric family of curves on the plane $\bR^2_{\bs \xi}$ and this family could serve as a local coordinate net, then in these coordinates the image of each integral manifold would be a rectangle. Therefore, the new coordinates satisfy the necessary condition of a separation and we obtain hope to build such a separation. At this moment it absolutely does not matter whether this coordinate net was built by a global change of variables or, as in the Kowalevski case, this net (the Zhukovsky net) depends in a very complicated way on the integral constants.

Studying the projections of $\mP$ onto the picture plane (or in more general case onto some picture manifold with dimension equal to the number of degrees of freedom) we must have in mind that $\mP$ itself may not have a global coordinate system. For example, in the rigid body dynamics $\mP=\bR^3{\times}SO(3)$ with the matrix group $SO(3)$ naturally embedded in $\bR^9$ by means of the so-called direction cosines. Any system of three coordinates either not covers $SO(3)$ or has inevitable singularities (like the Euler angles or their modifications). It is one of the reasons to consider $\mP$ as a subset in a vector space $V$ of sufficiently high dimension\footnote{The main reason is, of course, to have a possibility to introduce algebraic and rational functions.} presenting it as the level surface of some mapping
\begin{equation}
\notag \Gamma: V \to \bR^{\dim V-\dim\mP}.
\end{equation}
In the rigid body dynamics this mapping is formed by the so-called geometric integrals, and in the theory of the Euler equations on Lie algebras this role is given to the Casimir functions.

Combining (\ref{neq2_5}) with $\Gamma$ we obtain each integral manifold as a level surface of the new, enhanced, integral mapping $\mathcal{G}: V \to Z,$ where $V$ and $Z$ are real vector spaces and $\mF_{\mbf{f}}$ coincides with ${\mathcal{G}_{\mbf{z}} \dot = \mathcal{G}^{-1}(\mbf{z})}$, $\mbf{z} \in Z$. By an appropriate choice of coordinates we can represent ${V=X{\times}Y}$, where $X$ is the picture plane, and replace the projection of $\mP$ onto the picture plane with the projection of $V$ onto the first multiplier ${p_X: X{\times}Y \to X}$. This construction is easily generalized to the case when $X,Y,Z$ are smooth manifolds. In axially symmetric problems of the rigid body dynamics the manifold $X$ is often identified with the so-called Poisson sphere $S^2$. Then instead of $\A(\mbf{f})$ we have a similar region
\begin{equation}
\notag \A(\mbf{z}) = \operatorname{Im} p_X|_{\mathcal{G}_{\mbf{z}}}=p_X (\mathcal{G}^{-1}(\mbf{z})).
\end{equation}
To construct coordinate nets on the picture plane we need to investigate the boundaries of accessible regions, not only in the topological sense, but the curves that in geometry and in the wave fronts theory are called visible contours.

\begin{definition}\label{thdef23} {\rm \cite{Kh82}}. The generalized boundary of an accessible region $\mathcal{G}_\mbf{z}$ is the set of critical values of the mapping
\begin{equation}
\notag p_X|_{\mathcal{G}_\mbf{z}}: \mathcal{G}_\mbf{z} \to X.
\end{equation}
\end{definition}

To find generalized boundaries one needs, in particular, to solve the equation $\mathcal{G}=\mbf{z}$. The following result \cite{Kh82} helps to avoid this sometimes impossible procedure.

\begin{proposition}\label{th21} Let $\mbf{v} \in V$. Denote by $L(\mbf{v})$ the restriction of the operator
\begin{equation}
\notag T_\mbf{v} \mathcal{G} : X {\times} Y \to Z
\end{equation}
onto the subspace $Y$. The image of the point $\mbf{v} \in \mathcal{G}_\mbf{z}$ belongs to the generalized boundary of the accessible region $\A(\mbf{z})$ if and only if $\rk L(\mbf{v}) < \dim Z.$ \end{proposition}

Naturally, this statement is of practical use only if $\dim \mathcal{G}_\mbf{z} \geqslant \dim X$, which is equivalent to ${\dim Y \geqslant \dim Z}$. In global coordinates on $V$ and $Z$, the condition for $\rk L(\mbf{v})$ is easy to write in terms of the corresponding Jacobi matrix.

\subsection{Integrable problems of the rigid body dynamics} 
We tend to apply the above ideas to the system of equations describing the rotation of a rigid body about a fixed point in the linear potential force field
\begin{equation}
\label{neq2_12}
\begin{array}{l}
\displaystyle{\mbf{I}\frac{d{\bs \omega}}{dt}=\mbf{I}{\bs \omega}\times{\bs \omega}+\bf{r}_1 \times {\bs \alpha}+\mbf{r}_2 \times {\bs \beta},} \\
\displaystyle{\frac{d{\bs \alpha}}{dt}= {\bs \alpha}\times {\bs \omega},} \qquad \displaystyle{\frac{d{\bs \beta}}{dt}= {\bs \beta}\times {\bs \omega}.}
\end{array}
\end{equation}
Here ${\bs \omega}=(\omega_1,\omega_2,\omega_3)$, ${\bs \alpha}=(\alpha_1,\alpha_2,\alpha_3)$, ${\bs \beta}=(\beta_1,\beta_2,\beta_3)$ are the phase variables, the constant vectors $\mbf{r}_1,\mbf{r}_2 \in \bR^3$ and the diagonal constant matrix $\mbf{I}=\diag\{A,B,C\}$ are the physical parameters. The geometric integrals forming the corresponding mapping $\Gamma: \bR^9 \to \bR^3$ are ${\bs \alpha}^2, {\bs \beta}^2, {\bs \alpha}{\cdot} {\bs \beta}$. These are the Casimir functions for the Poisson brackets on $\bR^9$ bringing the equations (\ref{neq2_12}) to the Hamiltonian form \cite{BogEn}. If ${\bs \alpha}, {\bs \beta}$ (the intensity vectors of the forces) are linearly independent, the common level of the geometric integrals is a smooth 6-dimensional manifold $\mP$ endowed with a symplectic structure.

Put $\mbf{I}=\diag \{2,2,1\}$, $\mbf{r}_1 = (1,0,0)$, $\mbf{r}_2 = (0,1,0)$. Then equations (\ref{neq2_12}) describe the Kowalevski type top in a double force field. The first partial case of integrability (with restrictions on the phase variables) was found in \cite{BogRus1, BogRus2, BogEn}. In the special case of Yehia \cite{Yeh1} (with restrictions on the force fields $|{\bs \alpha}|=|{\bs \beta}|, {\bs \alpha}{\cdot} {\bs \beta}=0$) the system has a symmetry group and becomes reducible to a family of integrable systems with two degrees of freedom. The complete integrability of this problem in general was established in \cite{ReySem,ReySemRus}. More detailed presentation for the case of the Kowalevski gyrostat is given in \cite{Bob}.

\begin{definition} Let ${\bs \alpha},{\bs \beta}\in \bR^3$. The following constants
\begin{equation}\label{neq2_13}
p=\sqrt{{\bs \alpha}^2+{\bs \beta}^2} \geqslant 0 ,\qquad r=\sqrt[4]{({\bs \alpha}^2-{\bs \beta}^2)^2+4({\bs \alpha}{\cdot}{\bs \beta})^2}\geqslant 0
\end{equation}
are called the invariants of the pair $({\bs \alpha},{\bs \beta})$. \end{definition} Obviously, $p\geqslant r$ and $p=r$ if and only if ${\bs \alpha}{\times}{\bs \beta}=0$. It is the classical Kowalevski case. If $r=0$ we obtain the reducible Yehia case.

For $\theta \in \bR$ denote
\begin{equation}\label{neq2_14}
\Theta=\left(
\begin{array}{rr} \cos\theta & -\sin \theta \cr \sin\theta & \cos \theta \end{array}\right),\qquad \widetilde \Theta = \left(
\begin{array}{c|c} \Theta & 0 \\ \hline 0 & 1 \end{array}\right).
\end{equation}

\begin{proposition} \emph{\cite{Kh34}} The linear automorphism of $\bR^9$
\begin{equation}\label{neq2_15}
\begin{array}{l}
\begin{pmatrix}{\bs \alpha}'\cr {\bs \beta}'\end{pmatrix}=\Theta
\begin{pmatrix}{\bs \alpha}\cr {\bs \beta}\end{pmatrix}\widetilde \Theta^\top,\qquad  {\bs \omega}'= \widetilde \Theta \, {\bs \omega}. \end{array}
\end{equation}
preserves the equations $(\ref{neq2_12})$. An orbit of the restriction of this automorphism onto the space $\bR^6=\{({\bs \alpha},{\bs \beta})\}$ consists of all pairs having the same invariants. Each such orbit contains an orthogonal pair ${\bs \alpha} {\cdot} {\bs \beta}=0$. \end{proposition}

It follows that without loss of generality for irreducible problems we can consider the phase space $\mathcal{P} \subset \mbf{R}^9$ as given by the equations
\begin{equation}\label{neq2_16}
{\bs \alpha}^2=a^2,\qquad  {\bs \beta}^2=b^2, \qquad {\bs \alpha}{\cdot} {\bs \beta}=0\qquad (a>b>0).
\end{equation}
The proposition holds for an axially symmetric gyrostat in the same force field~\cite{KhND07}.

The integration of the whole generalized Kowalevski system has not been fulfilled. Separation of variables of the algebraic type is obtained only in the following three cases:

\begin{itemize} \item the subsystem with one degree of freedom consisting of the families of critical periodic motions in the first critical subsystem with two degrees of freedom; the subsystem itself was found by O.I.\,Bogoyavlensky in the works \cite{BogRus1, BogRus2}, the existence of periodic motions was proved in \cite{ZotRCD}, integration presented in \cite{Kh361};

\item the second critical subsystem with two degrees of freedom found in \cite{Odin}, separation of variables is pointed out in \cite{KhSavDan,KhSavMRC};

\item the third critical subsystem found in \cite{Kh34}, separation of variables is given in \cite{KhND06,Kh38,KhRCD09}. \end{itemize}

For the generalized Kowalevski top, let us present the system (\ref{neq2_12}) in the scalar form
\begin{equation}\label{neq2_17}
\begin{array}{c}
2\dot\omega _1   = \omega _2 \omega _3  + \beta _3 ,\;\quad 2\dot\omega _2 =  - \omega _1 \omega _3  - \alpha _3 ,\; \quad \dot\omega _3   = \alpha _2  - \beta _1 , \\
\dot\alpha _1   = \alpha _2 \omega _3  - \alpha _3 \omega _2 ,\; \quad \dot\beta _1   = \beta _2 \omega _3  - \beta _3 \omega _2, \\
\dot\alpha _2   = \alpha _3 \omega _1  - \alpha _1 \omega _3 ,\; \quad \dot\beta _2   = \beta _3 \omega_1  - \beta _1 \omega _3, \\
\dot\alpha_3   = \alpha_1 \omega_2  - \alpha_2 \omega_1 ,\; \quad \dot\beta_3   = \beta_1 \omega_2 - \beta_2 \omega_1. \end{array}
\end{equation}
Its first integrals are as follows \cite{BogRus1, ReySem}
\begin{equation}\label{neq2_18}
\begin{array}{l}
H = \omega _1^2  + \omega _2^2  + \ds{{\frac{1}{2}}}\omega _3^2 - (\alpha _1  + \beta _2 ), \\[2mm]
K = (\omega _1^2  - \omega _2^2  + \alpha _1 - \beta _2 )^2  + (2\omega _1 \omega _2  + \alpha _2  + \beta _1 )^2 , \\[2mm]
G = (\alpha _1 \omega _1  + \alpha _2 \omega _2  + \ds{{\frac{1}{2}}}\alpha _3 \omega _3 )^2  + (\beta _1 \omega _1  + \beta _2 \omega _2  + \ds{{\frac{1}{2}}}\beta _3 \omega _3 )^2  \\
\qquad {} + \omega _3 (\gamma _1 \omega _1  + \gamma _2 \omega _2 + \ds{{\frac{1}{2}}}\gamma _3 \omega _3 ) - \alpha _1 b^2  - \beta _2 a^2.
\end{array}
\end{equation}
Here $ \gamma_i$ stand for the components of the vector ${\bs{\alpha }}{\times}{\bs{\beta }}$. We denote the corresponding constants by $h,k,g$.

To make the description of critical subsystems compact, we use the complex change of variables \cite{Kh34} generalizing that of S.V.\,Kowalevski\footnote{These combinations of the real variables are in fact found in the Lax representation \cite{ReySem}.},
\begin{equation}\label{neq2_19}
\begin{array}{l}
x_1 = (\alpha_1  - \beta_2) + \ri (\alpha_2  + \beta_1),\quad x_2 = (\alpha_1  - \beta_2) - \ri (\alpha_2  + \beta_1 ), \\ y_1 = (\alpha_1  + \beta_2) + \ri (\alpha_2  - \beta_1), \quad y_2 = (\alpha_1  + \beta_2) - \ri (\alpha_2  - \beta_1), \\  z_1 = \alpha_3  + \ri \beta_3, \quad z_2 = \alpha_3  - \ri \beta_3,\\
w_1 = \omega_1  + \ri \omega_2 , \quad w_2 = \omega_1  - \ri \omega_2, \quad w_3 = \omega_3, \end{array}
\end{equation}
where $\ri^2=-1$. It brings the first integrals to the form
\begin{equation}\label{neq2_20}
\begin{array}{l}
\displaystyle{ H = \frac{1}{2}w_3^2 + w_1 w_2  - \frac{1}{2}(y_1 + y_2 ) },  \\[3mm]
\displaystyle{ K=(w_1^2 + x_1 )(w_2^2  + x_2 )},\\[3mm]
\displaystyle{ G = \frac{1}{4}(p^2  - x_1 x_2 )w_3^2 +\frac{1}{2}(x_2 z_1 w_1 +x_1 z_2 w_2 )w_3}  \\[2mm] \displaystyle{\qquad{} + \frac{1}{4}(x_2 w_1  + y_1 w_2 )(y_2 w_1 + x_1 w_2 ) - \frac{1}{4}p^2 (y_1  + y_2 )}\\[2mm]
\displaystyle{\qquad{} +\frac{1} {4}r^2 (x_1 + x_2 )}.
\end{array}
\end{equation}
Here according to (\ref{neq2_13}), (\ref{neq2_16}),
\begin{equation}
\label{neq2_21}
\begin{array}{c} p^2  = a^2  + b^2 ,\qquad r^2  = a^2  - b^2. \end{array}
\end{equation}
Equations (\ref{neq2_16}) take the form
\begin{eqnarray}
& z_1^2  + x_1 y_2  = r^2 ,\quad z_2^2  + x_2 y_1  = r^2 , \label{neq2_22}
\\ & x_1 x_2  + y_1 y_2  + 2z_1 z_2  = 2p^2 .\label{neq2_23}
\end{eqnarray}

The first critical subsystem $\mm \subset \mathcal{P}$ is defined by the equations
\begin{equation}
\notag w_1^2+x_1  = 0,\qquad w_2^2+x_2  = 0.
\end{equation}
On $\mm$ we have
\begin{equation}
\label{neq2_24} k = 0
\end{equation}
and we take $H$ and the Bogoyavlensky integral
\begin{equation}
\notag F = w_1 w_2 w_3+z_2 w_1+z_1 w_2
\end{equation}
as a pair of almost everywhere independent integrals.

The second critical subsystem $\mn \subset \mathcal{P}$ is defined by the equations
\begin{equation}
\label{neq2_25} x_1 x_2 w_3  - (x_2 z_1 w_1  + x_1 z_2 w_2 )  = 0, \quad \ds{ \frac{x_2}{x_1}(w_1^2+x_1)-\frac{x_1}{x_2}(w_2^2+x_2)}  = 0.
\end{equation}
On $\mn$, we choose the following pair of almost everywhere independent integrals:
\begin{equation}
\label{neq2_26}
\begin{array}{l} \displaystyle{M = \frac{1} {2 r^2}\bigl[\frac{x_2}{x_1}(w_1^2+x_1)+\frac{x_1}{x_2}(w_2^2+x_2)\bigr]},\\ \ds{L = \frac{1} {{\sqrt {x_1 x_2 } }}\bigl( w_1 w_2  + {{x_1 x_2  + z_1 z_2 }} M \bigr).} \end{array}
\end{equation}
The constants of the general integrals (\ref{neq2_18}) satisfy the equation
\begin{equation}
\label{neq2_27} (p^2h-2g^2)^2-r^4k=0.
\end{equation}

The third critical subsystem $\mo \subset \mathcal{P}$ is described as
\begin{equation}\label{neq2_28}
\begin{array}{l}
\displaystyle{\frac{w_2 x_1+w_1 y_2+w_3 z_1}{w_1}-\frac{w_1 x_2+w_2 y_1+w_3 z_2}{w_2} = 0,} \\
\displaystyle{(w_2 z_1+w_1 z_2)w_3^2+\Bigl[\frac{w_2 z_1^2}{w_1}+\frac{w_1 z_2^2}{w_2}+w_1 w_2(y_1+y_2)}\\
\qquad +\displaystyle{x_1 w_2^2+x_2 w_1^2\Bigr]w_3 +\frac{w_2^2 x_1 z_1}{w_1} + \frac{w_1^2 x_2 z_2}{w_2}}\\
\qquad + \displaystyle{ x_1 z_2 w_2+ x_2 z_1 w_1 +(w_1 z_2-w_2 z_1)(y_1-y_2)=0.}
\end{array}
\end{equation}
Here we take as a pair of almost everywhere independent integrals
\begin{equation}
\label{neq2_30}
\begin{array}{l}
\displaystyle{S=-\frac{1}{4} \bigl( \frac {y_2 w_1+x_1 w_2+z_1 w_3}{w_1}+\frac{x_2 w_1+y_1 w_2+z_2 w_3}{w_2} \bigr),} \\[3mm] \displaystyle{T =\frac{1}{2}[w_1(x_2 w_1+y_1 w_2+z_2 w_3)+w_2(y_2 w_1+x_1 w_2+z_1 w_3)]}\\
\phantom{T =} + x_1 x_2+z_1 z_2.
\end{array}
\end{equation}
The constants of (\ref{neq2_18}) on each manifold $\{{S=s},\,{T=\tau}\}$ satisfy the following equations
\begin{equation}
\label{neq2_31} h = \frac {p^2 -\tau}{2s}+ s, \quad k = \frac{\tau^2 - 2p^2\tau + r^4}{4s^2}+\tau, \quad g = \frac{p^4-r^4} {4s}+ \frac{1}{2}(p^2 -\tau)s.
\end{equation}

As $b$ tends to zero the critical subsystems turn into the famous classes of ``especially remarkable'' motions found by G.G.\,Appelrot \cite{Appel} in the Ko\-wa\-lev\-ski case. From the point of view of the phase topology of the integrable Hamiltonian system these critical subsystems form the critical set of the momentum mapping and equations (\ref{neq2_24}), (\ref{neq2_27}) and (\ref{neq2_31}) define the bifurcation diagram \cite{Kh34}.

\subsection{Two systems of local coordinates}
\subsubsection{The coordinates $s_1,s_2$}
We follow the advice of C.G.J.\,Jacobi \cite{Jacobi} to find first a convenient change of variables and then search for the problems in which it may be successfully applied.

\begin{figure}[ht]
\centering
\includegraphics[width=100mm, keepaspectratio]{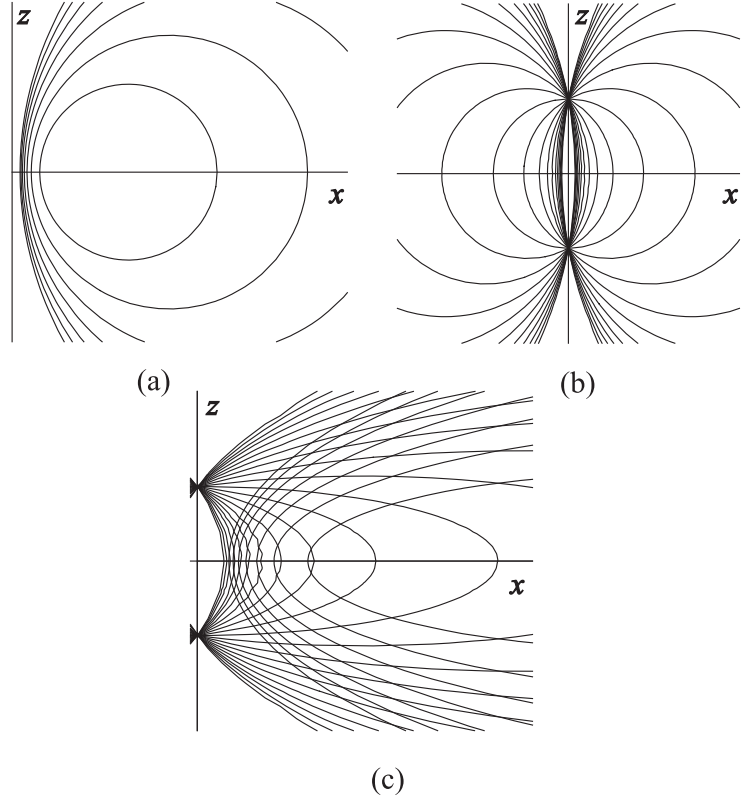} \caption{The coordinate net $(s_1,s_2)$.}\label{fig_s1s2} \end{figure}

Consider a number $r> 0$. On the plane $(x,z)$ introduce two one-pa\-ra\-met\-ric families of circles
\begin{equation}
\label{neq2_33} x^2  + z^2  + r^2  = 2s_1 x, \qquad x^2  + z^2  - r^2  = 2s_2 x.
\end{equation}
Except for the points of the axes $Ox,Oz$ these families provide the net of curvilinear coordinates $(s_1 ,s_2 )$ on the $(x,z)$-plane. In Fig.~\ref{fig_s1s2} we show: (a)~the curves $s_1  = \mathrm{const}$ for $s_1  \geqslant r$; (b)~the curves $s_2  = \mathrm{const}$ for $s_2$ in some interval about 0; (c)~joint net $(s_1 ,s_2)$ scaled $1{:}2$ for better view.
For the differentials we have
\begin{equation}
\label{neq2_34}
\begin{array}{l} \ds{ds_1  = \frac{{x^2  - z^2  - r^2 }}{{2x^2 }}dx + \frac{z}{x}dz}, \qquad \ds{ds_2  = \frac{{x^2  - z^2  + r^2 }}{{2x^2 }}dx + \frac{z}{x}dz}. \end{array}
\end{equation}

\begin{figure}[ht]
\centering
\includegraphics[width=60mm, keepaspectratio]{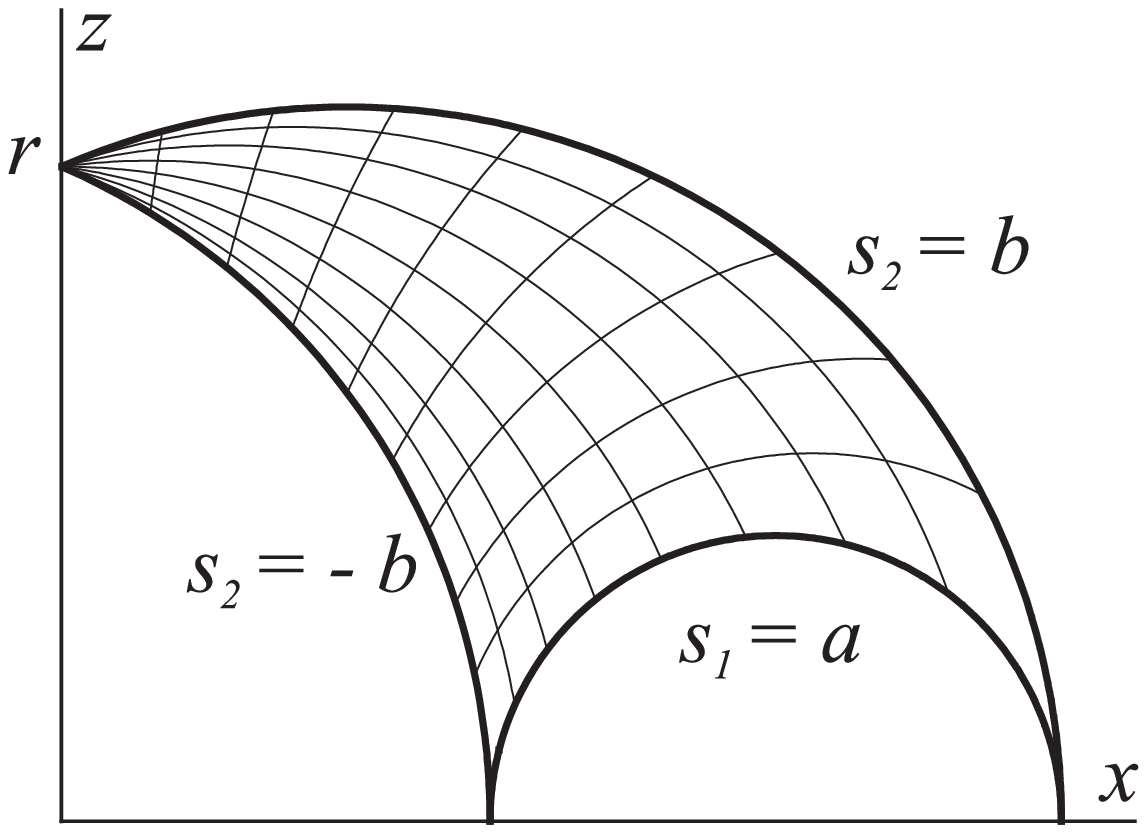}
\caption{The image of the phase space.}\label{fig_xznet}
\end{figure}

We now establish the connection with the rigid body configuration space (\ref{neq2_16}) and find the image of it on the $(s_1,s_2)$-plane. To this end, choose the parameters $p,r$ according to (\ref{neq2_21}). Introduce the variables $x,y,z$ as follows
\begin{equation}
\label{neq2_35} x^2  = x_1 x_2 ,\quad y^2  = y_1 y_2 ,\quad z^2  = z_1 z_2 .
\end{equation}
Then due to (\ref{neq2_23}) we obtain the ellipsoid
\begin{equation}
\notag x^2  + y^2  + 2z^2  = 2p^2.
\end{equation}
From (\ref{neq2_23}) express
\begin{equation}
\label{neq2_36} y_1 y_2  = 2p^2  - x^2  - 2z^2
\end{equation}
and represent (\ref{neq2_22}) in the form
\begin{equation}
\label{neq2_37}
\begin{array}{l} (z_1  + z_2 )^2  = 2r^2  - (x_1 y_2  + x_2 y_1 ) + 2z^2 ,  \\ (z_1  - z_2 )^2  = 2r^2  - (x_1 y_2  + x_2 y_1 ) - 2z^2 . \end{array}
\end{equation}
Elimination of the variables $z_1 ,z_2$ and the product $y_1 y_2$ in (\ref{neq2_22}), (\ref{neq2_23}), (\ref{neq2_36}) gives
\begin{equation}
\label{neq2_38} r^2 (x_1 y_2  + x_2 y_1 ) = r^4  + 2p^2 x^2  - (x^2  + z^2)^2 .
\end{equation}
Denote
$$
\Phi _ \pm  (x,z) = (x^2  + z^2  \pm r^2 )^2  - 2(p^2  \pm r^2 )x^2.
$$
Then (\ref{neq2_37}), (\ref{neq2_38}) yield the following expressions:
\begin{equation}
\label{neq2_39} r^2 (z_1  + z_2 )^2  = \Phi _ +  (x,z),\quad r^2 (z_1  - z_2 )^2  = \Phi _ -  (x,z).
\end{equation}
Therefore, the region for $x,z$ is defined by the inequalities
\begin{equation}
\label{neq2_40} \Phi _ +  (x,z) \geqslant 0,\quad \Phi _ -  (x,z) \leqslant 0.
\end{equation}
Define $s_1 ,s_2$ as in (\ref{neq2_33}) with the chosen $r$. Then
\begin{equation}
\label{neq2_41} s_1  = \ds{\frac{x^2  + z^2  + r^2} {2x}},\quad s_2  = \ds{\frac{x^2 + z^2 - r^2} {2x}}.
\end{equation}
The conditions (\ref{neq2_40}) take the form
\begin{equation}
\label{neq2_42} s_1^2  \geqslant a^2 ,\quad s_2^2  \leqslant b^2 .
\end{equation}
This rectangular area (in the above mentioned generalized sense) is the image of the phase space $\mP$ under the projection $({\bs \omega},{\bs \alpha},{\bs \beta}) \mapsto ({\bs \alpha},{\bs \beta}) \mapsto (s_1,s_2)$. The image of $\mP$ on the $(x,z)$-plane with the corresponding part of the $(s_1,s_2)$-net is shown in Fig.~\ref{fig_xznet} for the first quadrant.

\subsubsection{The coordinates $t_1,t_2$} 
\begin{figure}[ht]
\centering
\includegraphics[width=60mm, keepaspectratio]{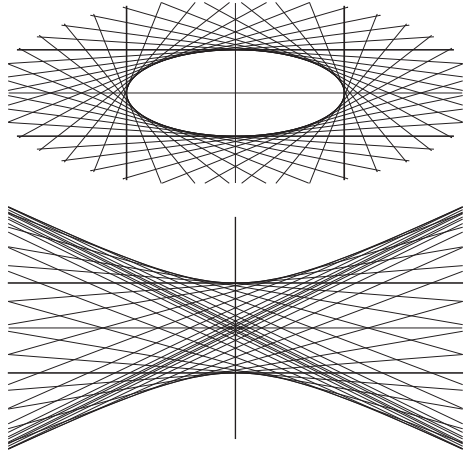}
\caption{The net of tangent lines to a second order curve.}\label{fig_t1t2}
\end{figure}

Consider the plane with coordinates $(\xi ,x)$ and the second order curve
\begin{equation}
\label{neq2_43} \frac{{\xi ^2 }}{\sigma } + \frac{{x ^2 }}{\tau } = 1,
\end{equation}
where $\tau ,\sigma$ are not simultaneously negative. Let $(\xi ,x )$ be an arbitrary point of the region
\begin{equation}
\label{neq2_44} \tau \xi ^2  + \sigma x ^2  > \tau \sigma.
\end{equation}
Then from this point we can draw exactly two tangent lines to the curve (\ref{neq2_43}). In other words, the system of tangent lines to the curve (\ref{neq2_43}) provides, in the region (\ref{neq2_44}), the net of straight lines such that each point belongs to exactly two lines of this net (see Fig.~\ref{fig_t1t2}). This fact is well-known in geometry.

To describe this net, introduce the variables $t_1 ,t_2 $ as the roots of the quadratic equation
\begin{equation}
\label{neq2_45} t^2  - \frac{{2\tau \xi }}{{\tau  - x ^2 }}t + \frac{{\tau \xi ^2 + \sigma x ^2 }}{{\tau  - x ^2 }} = 0.
\end{equation}
Its discriminant is positive in the region (\ref{neq2_44}), hence the following expressions are real,
\begin{equation}
\label{neq2_46} t_1  = \frac{{\tau \xi  + \mu x }}{{\tau  - x ^2 }},\quad t_2  = \frac{{\tau \xi  - \mu x }}{{\tau  - x ^2 }}.
\end{equation}
Here $\mu = \sqrt {\tau \xi ^2  + \sigma x ^2  - \tau \sigma }  \in {\mbf{R}}$.

Let us also introduce the connection between two coordinate systems $(t_1,t_2)$ and $(s_1,s_2)$. It will be useful below. Let $p,r,a,b$ satisfy (\ref{neq2_21}). In (\ref{neq2_43}) we put
\begin{equation}
\notag \sigma  = \tau ^2  - 2p^2 \tau  + r^4
\end{equation}
and let the coordinates $\xi ,x$ depend on $(s_1 ,s_2 )$
\begin{equation}
\label{neq2_47} \xi  = \frac{{s_1  + s_2 }}{{s_1  - s_2 }}r^2-\tau,\quad x  = \frac{{r^2 }}{{s_1  - s_2 }}.
\end{equation}
With respect to the variables (\ref{neq2_35}) we get
\begin{equation}
\label{neq2_48} \xi  = x^2  + z^2  - \tau.
\end{equation}
Then the boundary lines of the region (\ref{neq2_42})
\begin{equation}
\label{neq2_49} s_1  =  \pm a,\quad s_2  =  \pm b
\end{equation}
appear to be coordinate lines in the net $(t_1 ,t_2 )$:
\begin{equation}
\notag t_1  =  - \tau  \pm r^2 ,\quad t_2  =  - \tau  \mp r^2 .
\end{equation}

\section{Example of a geometry based elliptic separation}\label{sec3}
Here we present the result announced in \cite{KhSavDan,KhSavMRC}. Complete investigation including the study of the phase topology can be found in \cite{KhSav}.

Consider the second critical subsystem and fix the constants of the first integrals~(\ref{neq2_26}):
\begin{equation}
\label{neq3_1} M = m,\quad L = \ell .
\end{equation}
Eliminate from (\ref{neq2_25}), (\ref{neq2_22}), (\ref{neq2_23}) and (\ref{neq3_1}) the variables $w_1 ,w_2 ,w_3 ,y_1 ,y_2$ and the sum ${z_1  + z_2}$. The result is one equation containing $x_1 ,x_2$ and the product ${z_1 z_2}$:
\begin{equation}
\label{neq3_2} \ds{\frac{1}{x_1 x_2}}\Bigl[ m(x_1 x_2  + z_1 z_2) + \sqrt {r^4 m^2 - r^2 m (x_1+x_2) + x_1 x_2 }\Bigr]^2 = \ell^2.
\end{equation}
It defines a two-dimensional surface in $\bR^3(\operatorname{\mathfrak{Im}} x_1, \operatorname{\mathfrak{Re}} x_1, z)$. Introducing the polar angle $\phi$ by putting $x_1  + x_2  = 2x\cos \phi$ we obtain the equation on the cylinder $\bR^2(x,z){\times}S^1(\phi \mathop{\rm mod}\nolimits 2\pi)$:
\begin{equation}
\notag m(x^2  + z^2 ) + \sqrt {r^4 m^2  - 2r^2 m x\cos \phi  + x^2 }  \pm \ell x = 0.
\end{equation}
According to Proposition~\ref{th21}, the generalized boundary of the projection onto the $(x,z)$-plane of a set given by an equation of the type ${f(x,z,\phi) = 0}$ without any restrictions on $\phi $ is defined by the condition $\partial f/\partial \phi = 0$. In our case since the main equation contains only $\cos \phi$, this condition immediately gives $\sin \phi = 0$. Then on the generalized boundary
\begin{equation}
\notag m(x^2  + z^2 ) + |r^2 m \pm x| \pm \ell x = 0.
\end{equation}
Hence
\begin{equation}
\notag x^2  + z^2  \pm r^2  \pm \frac{{\ell  \pm 1}}{m}x = 0.
\end{equation}
Recalling the definition (\ref{neq2_33}) we see that for any combination of signs these curves on the $(x,z)$-plane become the coordinate lines of the $(s_1 ,s_2)$-net. In particular, all integral manifolds (\ref{neq3_1}) \textit{have rectangular images} on the $(s_1 ,s_2)$-plane (possibly unbounded if the value $x=0$ equivalent to $s_1=\infty$ is admissible). Let us show that this fact indeed leads to separation of variables.

Write the equations (\ref{neq2_17}) for the variables (\ref{neq2_19}):
\begin{equation}\label{neq3_3}
\begin{array}{ll}
{x'_1  =  - x_1 w_3  + z_1 w_1,} & {x'_2  = x_2 w_3  - z_2 w_2,} \\
{y'_1  =  - y_1 w_3  + z_2 w_1,} & {y'_2  = y_2 w_3  - z_1 w_2 ,}\\
{2z'_1  = x_1 w_2  - y_2 w_1,} & {2z'_2  =  - x_2 w_1 + y_1 w_2,}\\
2w'_1  =  - w_1 (w_3-\lambda) - z_1, & {}\\
2w'_2  = w_2 (w_3-\lambda) + z_2, & {}\\
2w'_3 = y_2  - y_1. & {}
\end{array}
\end{equation}
Here prime stands for $d/d(\ri\, t)$. Calculating the derivatives by virtue of these equations we get
\begin{equation}
\label{neq3_4}
\begin{array}{l} \ds{\bigl(\frac{1}{x}\bigr)' = \frac{1}{{2(x_1 x_2 )^{3/2} }}(x_1 z_2 w_2 - x_2 z_1 w_1 ),} \\[4mm] \ds{\bigl(\frac{{x^2  + z^2 }}{x}\bigr)' = \frac{1}{{2(x_1 x_2 )^{3/2} }}[x_1 z_1 (z_2^2  + y_2 x_1 )w_2  - x_2 z_2 (z_1^2  + y_1 x_2 )w_1 ].} \end{array}
\end{equation}
Taking into account (\ref{neq3_1}), we obtain from (\ref{neq2_34})
\begin{equation}
\label{neq3_5}
\begin{array}{l} \ds{s'_1  = \frac{{r^2 }}{{4x^2 }}(z_1  + z_2 )\bigl[\sqrt {\frac{{x_1 }}{{x_2 }}} w_2  - \sqrt {\frac{{x_2 }}{{x_1 }}} w_1 \bigr],} \\  \ds{s'_2  = \frac{{r^2 }}{{4x^2 }}(z_1  - z_2 )\bigl[\sqrt {\frac{{x_1 }}{{x_2 }}} w_2  + \sqrt {\frac{{x_2 }}{{x_1 }}} w_1 \bigr].}  \end{array}
\end{equation}
Denote
$$
R_1  = \sqrt {r^2 m - x_1 } ,\qquad R_2  = \sqrt {r^2 m - x_2 }
$$
and express from (\ref{neq2_25})
\begin{equation}
\label{neq3_6} \sqrt {\frac{{x_1 }}{{x_2 }}} w_2  = R_1 ,\quad \sqrt {\frac{{x_2 }}{{x_1 }}} w_1  = R_2 .
\end{equation}
Equation (\ref{neq3_2}) takes the form
\begin{equation}
\label{neq3_7} m(x^2+z^2) + R_1 R_2   = \ell x.
\end{equation}
It follows from the definition of $R_1 ,R_2$ that
\begin{equation}\label{neq3_8}
\begin{array}{l}
R_1^2  + R_2^2  = 2r^2 m - (x_1  + x_2 ), \\
R_1^2 R_2^2  = r^4 m^2 - 2r^2 m(x_1  + x_2 ) + x^2 .
\end{array}
\end{equation}
Therefore from (\ref{neq3_7}) we obtain
\begin{equation}\notag
\begin{array}{l}
R_1 R_2  = \ell x - m(x^2  + z^2 ), \\
\ds{R_1^2  + R_2^2  = \frac{1}{{r^2 m}}\bigl\{ [\ell x - m(x^2  + z^2 )]^2  - x^2 \bigr\}  + r^2 m.}
\end{array}
\end{equation}
Using the expressions of $x$, $x^2  + z^2$ via $s_1 ,s_2 $, we find
\begin{eqnarray}
& & \ds{(R_1  \pm R_2 )^2  = \frac{{r^2 }}{{m(s_1  - s_2 )^2 }}\{ (\ell ^2 - 1) - 2\ell m[(s_1  + s_2 ) }\nonumber \\
& & \qquad \mp (s_1  - s_2 )] + 2m^2 [(s_1^2  + s_2^2 ) \mp (s_1^2 - s_2^2 )]\} .\label{neq3_10}
\end{eqnarray}
Introduce the following notation:
\begin{equation}\label{neq3_11}
\begin{array}{l}
\ds{\Psi(s_1,s_2)=4m s_1 s_2-2\ell (s_1+s_2)+\ds{\frac{1}{m}}(\ell^2-1),}\\[3mm]
\ds{\Phi (s) = \Psi(s,s) = 4ms^2  - 4\ell s + \frac{1}{m}(\ell ^2  - 1).}
\end{array}
\end{equation}
Then from (\ref{neq3_10}),
\begin{equation}\label{neq3_12}
\begin{array}{l}
\ds{    R_1  + R_2  = \frac{r}{{s_1  - s_2 }}\sqrt {\Phi (s_2 )},} \qquad  \ds{R_1  - R_2  = \frac{r}{{s_1  - s_2 }}\sqrt {\Phi (s_1 )}}
\end{array}
\end{equation}
(the radicals may have arbitrary signs). From (\ref{neq2_39}) we get
\begin{equation}\notag
\begin{array}{l}
\ds{(z_1  + z_2 )^2  = \frac{{4r^2 }}{{(s_1  - s_2 )^2 }}(s_1^2  - a^2 ),}\qquad
\ds{(z_1  - z_2 )^2  = \frac{{4r^2 }}{{(s_1  - s_2 )^2 }}(s_2^2  - b^2 ).}
\end{array}
\end{equation}
Hence
\begin{equation}
\label{neq3_14} \ds{z_1  + z_2  = \frac{{2r}}{{s_1  - s_2 }}\sqrt {s_1^2  - a^2 },} \qquad \ds{z_1  - z_2  = \frac{{2r}}{{s_1  - s_2 }}\sqrt {s_2^2  - b^2 },}
\end{equation}
where the signs of the radicals, similar to (\ref{neq3_12}), are arbitrary.

Substitute (\ref{neq3_6}), (\ref{neq3_11}), (\ref{neq3_14}) in (\ref{neq3_5}) and return to the real time derivative to obtain the equations describing the time evolution of the auxiliary variables $s_1 ,s_2 $:
\begin{equation}
\label{ds1ds2}
\begin{array}{l} \ds{\frac{{ds_1 }}{{dt}} = \frac{1}{2}\sqrt {\mstrut s_1^2  - a^2 } \sqrt {\mstrut  - \Phi (s_1 )} ,} \qquad  \ds{\frac{{ds_2 }}{{dt}} = \frac{1}{2}\sqrt {\mstrut b^2  - s_2^2 } \sqrt {\mstrut \Phi (s_2 )}.} \end{array}
\end{equation}
In these equations the variables are separated and the equations themselves are easily integrated in Jacobi's elliptic functions.

To accomplish the separation, it is necessary to present the expressions of all phase variables in terms of the separated ones. The values $z_1,z_2$ are found from (\ref{neq3_14}). From (\ref{neq2_33})
\begin{equation}
\notag
\begin{array}{l} x=\ds{\frac{r^2}{s_1-s_2}}, \qquad x^2+z^2=r^2 \ds{\frac{s_1+s_2}{s_1-s_2}}, \end{array}
\end{equation}
then $x_1,x_2$ are found from (\ref{neq3_8}) and (\ref{neq3_12}). The equations (\ref{neq2_22}) serve to define $y_1,y_2$, and $w_1,w_2$ are expressed from (\ref{neq3_6}). Finally, the variable $w_3$ is defined from the first equation (\ref{neq2_25}). Thus, for the complex phase variables we have
\begin{equation}\notag
\begin{array}{l}
\ds{x_1  =  - \frac{{r^2 }}{{2(s_1  - s_2 )^2 }}[\Psi (s_1 ,s_2 ) + \sqrt {\Phi (s_1 )\Phi (s_2 )} ], }\\[3mm]
\ds{x_2  =  - \frac{{r^2 }}{{2(s_1  - s_2 )^2 }}[\Psi (s_1 ,s_2 ) - \sqrt {\Phi (s_1 )\Phi (s_2 )} ], }\\[3mm]
\ds{y_1  = 2\frac{{(2s_1 s_2  - p^2 ) - 2\sqrt {(s_1^2  - a^2 )(s_2^2  - b^2 )} }}{{\Psi (s_1 ,s_2 ) - \sqrt {\Phi (s_1 )\Phi (s_2 )} }}, }\\[4mm]
\ds{y_2  = 2\frac{{(2s_1 s_2  - p^2 ) + 2\sqrt {(s_1^2  - a^2 )(s_2^2  - b^2 )} }}{{\Psi (s_1 ,s_2 ) + \sqrt {\Phi (s_1 )\Phi (s_2 )} }}, }\\[4mm]
\ds{z_1  = \frac{r}{{s_1  - s_2 }}(\sqrt {s_1^2  - a^2 }  + \sqrt {s_2^2  - b^2 } ), }\\[3mm]
\ds{z_2  = \frac{r}{{s_1  - s_2 }}(\sqrt {s_1^2  - a^2 }  - \sqrt {s_2^2  - b^2 } ) }
\end{array}
\end{equation}
and
\begin{equation}\notag
\begin{array}{l}
\ds{w_1  = r\frac{{\sqrt {\Phi (s_2 )}  - \sqrt {\Phi (s_1 )} }}{{\Psi (s_1 ,s_2 ) - \sqrt {\Phi (s_1 )\Phi (s_2 )} }}, }\\[4mm]
\ds{w_2  = r\frac{{\sqrt {\Phi (s_2 )}  + \sqrt {\Phi (s_1 )} }}{{\Psi (s_1 ,s_2 ) + \sqrt {\Phi (s_1 )\Phi (s_2 )} }}, }\\[4mm]
\ds{w_3  = \frac{1}{{s_1  - s_2 }}[\sqrt {(s_2^2  - b^2 )\Phi (s_1 )}  - \sqrt {(s_1^2  - a^2 )\Phi (s_2 )} ].}
\end{array}
\end{equation}
Denote
\begin{gather}\notag
\begin{array}{l} \displaystyle{S_1 = \sqrt {s_1^2 - a^2},}\quad \displaystyle{\varphi_1 = \sqrt {\mstrut - \Phi (s_1)},}\quad \displaystyle{S_2 = \sqrt {b^2 - s_2^2},}\quad \displaystyle{\varphi_2 = \sqrt {\mstrut\Phi (s_2)}.} \end{array} \end{gather} Then
\begin{gather}\label{neq3_15}
\begin{array}{l}
\displaystyle{\alpha _1  = \frac{\mathstrut 1} {{2(s_1  - s_2 )^2 }}[(s_1 s_2 - a^2 )\Psi + S_1 S_2 \varphi _1 \varphi _2 ], }\\
\displaystyle{\alpha _2  = \frac{\mathstrut 1} {{2(s_1  - s_2 )^2 }}[(s_1 s_2 - a^2)\varphi _1 \varphi _2  - S_1 S_2 \Psi ], }\\
\displaystyle{\beta _1  =  - \frac{\mathstrut 1} {{2(s_1  - s_2 )^2 }}[(s_1 s_2  - b^2)\varphi _1 \varphi _2  - S_1 S_2 \Psi ], }\\
\displaystyle{\beta _2  = \frac{\mathstrut 1} {{2(s_1  - s_2 )^2 }}[(s_1 s_2 - b^2 )\Psi + S_1 S_2 \varphi _1 \varphi _2 ], }\\
\displaystyle{\alpha _3  = \frac{\mathstrut r} {{s_1  - s_2 }}S_1 ,\quad \beta _3  =\frac{r} {{s_1  - s_2 }}S_2, }\\ \displaystyle{\omega _1  = \frac{\mathstrut r} {{2(s_1  - s_2 )}}(\ell - 2ms_1 )\varphi _2,\quad \omega _2  = \frac{r} {{2(s_1 - s_2 )}}(\ell  - 2ms_2)\varphi _1, }\\
\displaystyle{\omega _3  = \frac{\mathstrut 1} {{s_1  - s_2 }}(S_2 \varphi _1 - S_1 \varphi _2 )}.
\end{array}
\end{gather}
Together with the differential equations (\ref{ds1ds2}) this gives the complete analytical solution of the problem.

\section{Example of a hyperelliptic separation}\label{sec4}
\subsection{The geometry of accessible regions} 
Let us consider the dynamical system defined by the equations (\ref{neq2_17}) on the invariant almost everywhere four-dimensional manifold~$\mo$ cut in $\bR^9$ by five equations (\ref{neq2_22}), (\ref{neq2_23}), (\ref{neq2_28}). This system is Hamiltonian with two degrees of freedom except for the zero-measure set of points where either the equations (\ref{neq2_28}) fail to be independent, or the induced symplectic structure degenerates~\cite{KhND06}. The possibility of separation of variables was predicted in \cite{Kh35} basing on the above shown geometrical speculations, but it took a long time to obtain the final analytical expressions, complex first \cite{KhND06} and finally real \cite{Kh38}.

The functions $S,T$ defined by (\ref{neq2_30}) provide the complete set of integrals in involution on $\mo$. In particular, the equations of the integral manifold
\begin{equation}
\notag \{ \zeta \in \mP :H(\zeta ) = h,K(\zeta ) = k,G(\zeta ) = g\}
\end{equation}
with the functions (\ref{neq2_20}) are replaced with the invariant relations (\ref{neq2_28}) and the system of equations
\begin{equation}
\label{neq4_1} S = s,\; T = \tau.
\end{equation}

Using the variables (\ref{neq2_19}) we transform the system (\ref{neq2_28}), (\ref{neq4_1}) to the form
\begin{equation}\label{neq4_2}
\begin{array}{l}
(y_2  + 2s)w_1  + x_1 w_2  + z_1 w_3  = 0,  \\
x_2 w_1  + (y_1  + 2s)w_2  + z_2 w_3  = 0,  \\
x_2 z_1 w_1  + x_1 z_2 w_2  + (\tau  - x_1 x_2 )w_3  = 0,  \\
2sw_1 w_2  - (x_1 x_2  + z_1 z_2 ) + \tau  = 0.
\end{array}
\end{equation}
Using the variables (\ref{neq2_19}) with the corresponding dynamics (\ref{neq3_3}) note that despite the complexity of these variables their space is 9-dimensional since the pairs ${(x_1,x_2)}$, ${(y_1,y_2)}$ and ${(z_1,z_2)}$ are complex conjugate and $w_3$ is in fact real. Seven relations (\ref{neq4_2}), (\ref{neq2_22}) and (\ref{neq2_23}) define in this space the integral manifold $\mathcal{F}_{s,\tau}$. If the latter set does not have any points of dependence of the functions $S,T$, then it will consist of two-dimensional tori bearing the quasi-periodic motions.

From what follows, we exclude the four points of equilibria on $\mP$. At non-trivial motions we have ${\bs \omega} \ne 0$ for almost all time moments. Then the determinant of the first three equations (\ref{neq4_2}) in $w_i$ ${(i=1,2,3)}$ is identically zero. Eliminate in this condition the variables $z_1^2, z_2^2$ and the product $y_1 y_2$ with (\ref{neq2_22}), (\ref{neq2_23}) and (\ref{neq2_36}) to obtain
\begin{equation}\label{neq4_3}
\begin{array}{l}
2s[(r^2 x_1  - \tau y_1 ) + (r^2 x_2  - \tau y_2 )] = - r^2 (x_1 y_2  + x_2 y_1 ) \\
\qquad \qquad + 2[2s^2 (\tau  - x^2 ) + p^2 (\tau  + x^2 ) - \tau (x^2  + z^2 )].
\end{array}
\end{equation}
On the other hand, it follows  from (\ref{neq2_35}), (\ref{neq2_36}) that
\begin{equation}
\label{neq4_4} (r^2 x_1  - \tau y_1 )(r^2 x_2  - \tau y_2 ) = r^4 x^2  + \tau (2p^2 - x^2  - 2z^2 ) - r^2 \tau (x_1 y_2  + x_2 y_1 ).
\end{equation}
Denote
\begin{equation}
\label{neq4_5} \sigma  = \tau ^2  - 2p^2 \tau  + r^4 ,\quad \chi  = \sqrt k \geqslant 0.
\end{equation}
We have,  from the second equation (\ref{neq2_31}),
\begin{equation}
\notag 4s^2 \chi ^2  = \sigma  + 4s^2 \tau .
\end{equation}
Introduce the complex conjugate pair
\begin{equation}
\label{neq4_6} \mu _1  = r^2 x_1  - \tau y_1 ,\quad \mu _2  = r^2 x_2  - \tau y_2
\end{equation}
and denote for brevity
\begin{equation}
\label{neq4_7} \xi= x^2+z^2-\tau.
\end{equation}
Eliminate in (\ref{neq4_3}), (\ref{neq4_4}) the term $x_1 y_2 + x_2 y_1$ by (\ref{neq2_38}):
\begin{equation}
\label{neq4_8}
\begin{array}{c} 2s(\mu _1  + \mu _2 ) = \xi^2  - 4s^2 (x^2  - \tau) - \sigma, \quad  \mu _1 \mu _2  = \tau \xi^2  + \sigma x^2  - \tau \sigma . \end{array}
\end{equation}
These two equations in four variables describe the integral manifold $\mathcal{F}_{s,\tau}$ in the space $(x,\xi,\mu_1,\mu_2)$ ($\mu_2=\overline{\mu_1}$). According to Proposition~\ref{th21}, at the points covering the generalized boundary of the projection of $\mathcal{F}_{s,\tau}$ onto the $(x,\xi)$-plane the Jacobian of the system with respect to the variables $\mu_1,\mu_2$ equals zero, or in view of (\ref{neq4_5}),
$$
(\xi -2sx-2s\chi)(\xi -2 s x+2 s\chi)(\xi +2sx-2s\chi)(\xi +2sx+2s\chi)=0.
$$
So, the straight lines
\begin{equation}\label{neq4_10}
\xi \pm 2sx \pm 2s\chi=0
\end{equation}
together with (\ref{neq2_49}) bound the projection of the integral manifold onto the $(x,\xi)$-plane.
Comparing (\ref{neq4_7}) with (\ref{neq2_48}) recall that the boundary lines of the region (\ref{neq2_42}) after the change (\ref{neq2_41}) become tangent to the second order curve
\begin{equation}\label{neq4_9}
\tau \xi^2  + \sigma x^2  - \tau \sigma =0.
\end{equation}
The resultant in $\xi$ of any of the equation (\ref{neq4_10}) with (\ref{neq4_9}) depending on the chosen signs evaluates to $4s^2(\tau \pm \chi x)^2$. Hence all straight lines (\ref{neq4_10}) are also tangent to the curve (\ref{neq4_9}), namely, at the points with $x=\mp \tau/\chi$.

Choose $\lambda _1  = \sqrt {2s\mu _1  + \sigma }$, $\lambda _2 = \sqrt {2s\mu _2  + \sigma }$ to be complex conjugate. Then the system (\ref{neq4_8}) simplifies $$ (\lambda _1  + \lambda _2 )^2  = \Xi _ +  (x,z),\quad (\lambda _1 - \lambda _2 )^2  = \Xi _ -  (x,z). $$ Here $$ \Xi _ \pm  (x,z) = (x^2  + z^2  - \tau  \pm 2s\chi )^2  - 4s^2 x^2, $$ and the solution exists if and only if
\begin{equation}
\label{neq4_11} \Xi _ +  (x,z) \geqslant 0,\quad \Xi _ -  (x,z) \leqslant 0.
\end{equation}
The inequalities (\ref{neq2_40}), (\ref{neq4_11}) describe the corresponding accessible region $\A(s,\tau)$ on the $(x,z)$-plane. We now present this region in the plane of the variables~$(s_1,s_2)$. Express
\begin{equation}\label{neq4_12}
\begin{array}{c}
\displaystyle{x^2  + z^2  - \tau  = [s_1  + s_2  - \ds{\frac{\tau}{r^2 }} (s_1  - s_2 )]x,}  \\
\displaystyle{\Xi _ +  (x,z) = x^2 \Lambda _ +  \Lambda _ -  ,\quad \Xi _ -  (x,z) = x^2 {\rm M}_ +  {\rm M}_-,}
\end{array}
\end{equation}
where
$$
\begin{array}{l}
\displaystyle{\Lambda _ \pm  (s_1 ,s_2 ) = s_1  + s_2  - \ds{\frac{\tau  - 2s\chi}{r^2 }} (s_1  - s_2 ) \pm 2s, } \\[4mm] \displaystyle{{\rm M}_ \pm  (s_1 ,s_2 ) = s_1  + s_2  - \ds{\frac{\tau  + 2s\chi}{r^2 }}(s_1  - s_2 ) \pm 2s. }
\end{array}
$$
From (\ref{neq4_11}), (\ref{neq4_12}) we have
\begin{equation}
\label{neq4_13} \Lambda _ +  (s_1 ,s_2 )\Lambda _ -  (s_1 ,s_2 ) \geqslant 0,\quad {\rm M}_ +  (s_1 ,s_2 ){\rm M}_ -  (s_1 ,s_2 ) \leqslant 0.
\end{equation}

\begin{figure}[ht]
\centering \includegraphics[width=100mm, keepaspectratio]{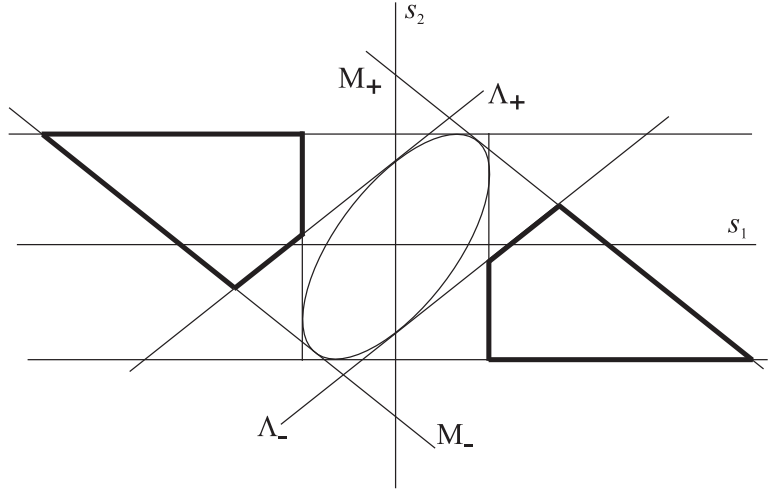}
\caption{Example of an accessible region.}\label{fig_accr} \end{figure}

The lines $\Lambda _ \pm = 0,\;{\rm M}_ \pm = 0$ form a parallelogram. The solutions of (\ref{neq4_13}) fill two half-strips attached to its sides belonging to the lines $\Lambda _ \pm   = 0$. Moreover, the solutions of (\ref{neq2_42}) fill two horizontal half-strips attached to the sides $s_2=\pm b$ of the rectangle with vertices $s_1 = \pm a$, $s_2 = \pm b$. The intersection of these sets give the accessible region for the variables~$s_1,s_2$ (e.g. see Fig.~\ref{fig_accr} for $a = 1$, $b = 0.4$, $\tau = 1.2$, $s =  - 0.6$). With respect to these variables all boundary lines are also tangent to the second order curve, which is the image of the curve (\ref{neq4_9}) under the change (\ref{neq2_47}).

The presented geometry on the planes $(x,z)$, $(x,\xi)$ and $(s_1,s_2)$ shows that it is possible to separate the variables by use of the second coordinate system (\ref{neq2_46}), where the notation have the sense given in this section.

\subsection{Parametric equations of integral manifolds}
At the two-dimensional surface (\ref{neq4_1}) the phase variables (\ref{neq2_19}) can be expressed in terms of two. At this moment it is convenient to choose for that the variables~$x,\xi$. Then from the equations (\ref{neq4_8}) we have
\begin{equation}\notag
\begin{array}{l}
\displaystyle{(\mu_1- 2 s \tau) + (\mu_2 - 2 s \tau) = \frac{\xi^2}{2s} - 2 s (x^2+\chi^2)},\\[2.5mm]
\displaystyle{(\mu_1- 2 s \tau)(\mu_2 - 2 s \tau)=4 s^2 \chi^2 x^2}.
\end{array}
\end{equation}
Hence
\begin{equation}
\label{neq4_14} \displaystyle{\mu_1=2s\tau+\frac{1}{8s}(\sqrt{\Psi_1}+\sqrt{\Psi_2})^2, \quad \mu_2=2s\tau+\frac{1}{8s}(\sqrt{\Psi_1}-\sqrt{\Psi_2})^2.}
\end{equation}
Here we denote
\begin{equation}\notag
\Psi_1(x,\xi) = \xi^2  -4s^2(x + \chi )^2, \quad \Psi_2(x,\xi) = \xi^2  -4s^2(x -\chi )^2.
\end{equation}
From (\ref{neq2_36}) we have
\begin{equation}
\notag y_1 y_2=2p^2-2(\xi+\tau)+x^2.
\end{equation}
Consider this equation together with the definition $x_1 x_2=x^2$ and the equations (\ref{neq4_6}). In view of (\ref{neq4_14}) we get
\begin{eqnarray}
& & \begin{array}{l}
    \displaystyle{x_1=\frac{2s}{r^2} \frac{4r^4(x^2-\tau)+\tau(\sqrt{\mathstrut \Phi_1}+\sqrt{\mathstrut \Phi_2})^2} {16s^2\tau+(\sqrt{\mathstrut \Psi_1}-\sqrt{\mathstrut \Psi_2})^2},} \\[4mm]
    \displaystyle{x_2=\frac{2s}{r^2} \frac{4r^4(x^2-\tau)+\tau(\sqrt{\mathstrut \Phi_1}-\sqrt{\mathstrut \Phi_2})^2} {16s^2\tau+(\sqrt{\mathstrut \Psi_1}+\sqrt{\mathstrut \Psi_2})^2},}
    \end{array} \label{neq4_15} \\
& & \begin{array}{l}
     \displaystyle{y_1= 2s \frac{4[2 \tau \xi - \tau (x^2-\tau)+ \sigma] - (\sqrt{\mathstrut \Phi_1}-\sqrt{\mathstrut \Phi_2})^2} {16s^2\tau+(\sqrt{\mathstrut \Psi_1}-\sqrt{\mathstrut \Psi_2})^2},} \\[4mm]
     \displaystyle{y_2=2s \frac{4[2 \tau \xi - \tau (x^2-\tau)+ \sigma] - (\sqrt{\mathstrut \Phi_1}+\sqrt{\mathstrut \Phi_2})^2} {16s^2\tau+(\sqrt{\mathstrut \Psi_1}+\sqrt{\mathstrut \Psi_2})^2}.}
     \end{array} \label{neq4_16}
\end{eqnarray}
Here the polynomials
\begin{equation}\label{neq4_17}
\begin{array}{c}
\Phi_1(x,\xi) = (\xi +\tau + r^2 )^2  - 2(p^2  + r^2 )x^2, \\
\Phi_2(x,\xi) = (\xi +\tau  - r^2 )^2  - 2(p^2  - r^2 )x^2
\end{array}
\end{equation}
are obtained from $\Phi_\pm$ after substituting $z^2$ with its value from (\ref{neq4_7}). In the sequel such substitution is supposed by default. From (\ref{neq2_37}) we then get
\begin{equation}
\label{neq4_18} \displaystyle{z_1=\frac{1}{2r} (\sqrt{\mathstrut \Phi_1} + \sqrt{\mathstrut \Phi_2}), \quad z_2=\frac{1}{2r} ( \sqrt{\mathstrut \Phi_1} - \sqrt{\mathstrut \Phi_2}).}
\end{equation}
As a result we have the expressions for all phase variables from the so-called configuration group, i.e., determining the orientation of the body. The radicals in (\ref{neq4_15}), (\ref{neq4_16}), (\ref{neq4_18}) have arbitrary signs.

To find the variables $w_i$, which are in charge of the angular velocity components, we use the expressions (\ref{neq2_20}) for the general integrals $K,H$ and the relations (\ref{neq2_31}), which define the constants $k,h$ in terms of $s,\tau$. In view of (\ref{neq4_5}) we have
\begin{eqnarray}
& & K=(w_1^2+x_1)(w_2^2+x_2)=\chi^2, \label{neq4_19}\\
& & H=\displaystyle{\frac{1}{2}w_3^2+w_1 w_2-\frac{1}{2}(y_1+y_2)=s+\frac{p^2-\tau}{2s}}. \label{neq4_20}
\end{eqnarray}
From (\ref{neq4_19}), (\ref{neq2_30}) we obtain the system
\begin{equation}\notag
\begin{array}{l}
\displaystyle{x_2 w_1^2+x_1 w_2^2=-\frac{1}{4s^2}[\xi^2+4s^2(\chi^2+x^2)],}\qquad \displaystyle{ w_1 w_2 = \frac{\xi}{2s}}. \end{array}
\end{equation}
Hence
\begin{equation}
 \displaystyle{w_1 = \frac{i}{4s\sqrt{x_2}}(\sqrt{\mathstrut \Theta_1}+\sqrt{\mathstrut \Theta_2}), \quad w_2 = \frac{i}{4s\sqrt{x_1}}(\sqrt{\mathstrut \Theta_1}-\sqrt{\mathstrut \Theta_2}),} \label{neq4_21}
\end{equation}
where
\begin{equation}
 \Theta_1(x,\xi) =(\xi-2s x)^2-4 s^2 \chi ^2,\quad \Theta_2(x,\xi) =(\xi + 2s x)^2-4 s^2 \chi ^2. \label{neq4_22}
\end{equation}
Substitute (\ref{neq4_21}), (\ref{neq4_16}) into (\ref{neq4_20}) to get
\begin{equation}\label{neq4_23}
\displaystyle{w_3^2  = \frac{1}{4s\mu ^2 }\bigl(P-\sqrt {\Phi _1 \Phi _2 \Psi _1 \Psi _2}\bigr),}
\end{equation}
where
\begin{equation}\notag
\begin{array}{l}
P =4s^2 (x^2  - \chi ^2 )[2(\tau  - p^2 )x^2  - \tau ^2  + r^4 ]+8s^2 [(\tau  - 2\chi ^2 )x^2  + \tau \chi ^2 ]\xi \\ \phantom{ P = }- 2[(\tau  - p^2  - 2s^2 )x^2  + \tau (p^2  - 2s^2 ) - r^4 ]\xi ^2- 2\tau \xi ^3- \xi ^4.
\end{array}
\end{equation}
Denote
\begin{equation}
\notag Q = (\xi  + \tau  + 2s^2  - p^2 )^2  - 4s^2 x^2  - (p^2  - 2s^2 )^2 + r^4.
\end{equation}
The following identity holds
$$
P^2 -\Phi _1 \Phi _2 \Psi _1 \Psi _2 = 4x^2 (\tau \xi ^2  + \sigma x^2  - \tau \sigma )Q^2.
$$
Using the second equation~(\ref{neq4_8}) we find from (\ref{neq4_23})
\begin{equation}
\label{neq4_24} \displaystyle{w_3  = \frac {1} {2\sqrt {2s} \mu }(\sqrt {P_1 } - \sqrt {P_2 } ),}
\end{equation}
where
\begin{equation}
\label{neq4_25} P_1  = P + 2x\mu Q,\quad P_2  = P - 2x\mu \delta, \qquad \mu^2=\mu_1 \mu_2.
\end{equation}
Simultaneously, from (\ref{neq4_8}),
\begin{equation}
\label{neq4_26} \mu^2=\tau \xi ^2  + \sigma x^2  - \tau \sigma.
\end{equation}
Thus, the equations (\ref{neq4_15}), (\ref{neq4_16}), (\ref{neq4_18}), (\ref{neq4_21}) and (\ref{neq4_24}) give the dependencies of all phase variables (\ref{neq2_19}) on the two variables $x,\xi$ taken as the intermediate ones.

\subsection{The change of variables}
Introduce the local coordinates $t_1,t_2$ according to (\ref{neq2_46}) on the $(x,\xi)$-plane. The discriminant of the corresponding equation (\ref{neq2_45}) now coincides with the value $\mu^2$ in (\ref{neq4_25}) and is, therefore, non-negative. So the variables $t_1,t_2$ are real. Solve (\ref{neq4_26}), (\ref{neq2_46}) with respect to $x, \xi, \mu$,
\begin{equation}
\label{neq4_27}
\begin{array}{c} \displaystyle{x = \frac{\sqrt{\tau}(U_1+U_2)}{t_1+t_2},} \quad \displaystyle{\xi = \frac{t_1 t_2+\sigma - U_1 U_2}{t_1+t_2},} \quad \displaystyle{\mu=\frac{\sqrt{\tau}(t_2 U_1 - t_1 U_2)}{t_1+t_2}.} \end{array}
\end{equation}
Here we denote for brevity
\begin{equation}
\label{neq4_28} U_1=\sqrt{t_1^2-\sigma}, \quad U_2=\sqrt{t_2^2-\sigma}.
\end{equation}
Note that the signs of these radicals should be chosen arbitrary to supply, at given $t_1, t_2$, all possible triples $x, \xi, \mu$ in (\ref{neq4_27}) satisfying the system (\ref{neq4_26}), (\ref{neq2_46}). The following easily checked identities simplify some calculation below:
\begin{eqnarray}
& & \displaystyle{x \mu =\frac{t_1-t_2}{t_1+t_2}\,\tau \xi,} \label{neq4_29} \\
& & \displaystyle{\frac{(t_1 t_2+\sigma + U_1 U_2)(t_1 t_2+\sigma - U_1 U_2)} {(t_1+t_2)^2} = \sigma,} \label{neq4_30} \\
& & \displaystyle{\frac{t_1 t_2+\sigma - U_1 U_2} {t_1 t_2-\sigma - U_1 U_2} = \frac{(U_1+U_2)^2}{(t_1-t_2)^2}.} \label{neq4_31}
\end{eqnarray}

We begin with $\mu_1, \mu_2$ defining the denominators in~(\ref{neq4_15}), (\ref{neq4_16}). Put $\psi = 16 s^2 \tau + \Psi_1 + \Psi_2$. According to (\ref{neq4_14}) we have
\begin{equation}\notag
\begin{array}{c}
\displaystyle{\mu_1=\frac{1}{8s}( \psi-2 \sqrt{\Psi_1 \Psi_2}), \qquad \mu_2=\frac{1}{8s}(\psi+2 \sqrt{\Psi_1 \Psi_2})}. \end{array}
\end{equation}
By definition $\psi^2 - 4 \Psi_1 \Psi_2 = 64 s^2 \mu^2$, and we write
$$
2 \sqrt{\mstrut \Psi_1 \Psi_2}=\sqrt{\mstrut  \psi + 8 s \mu} \sqrt{\mstrut  \psi - 8 s \mu}.
$$
Hence
\begin{equation}\notag
\begin{array}{l}
\displaystyle{\mu_1=\frac{1}{16 s}(\sqrt{\psi + 8 s \mu}-\sqrt{\psi - 8 s \mu})^2,} \\[3mm]
\displaystyle{\mu_2=\frac{1}{16 s}(\sqrt{\psi + 8 s \mu}+\sqrt{\psi - 8 s \mu})^2.}
\end{array}
\end{equation}
Substituting (\ref{neq4_27}) we get
\begin{equation}\notag
\begin{array}{l}
\displaystyle{\psi+8s\mu=\frac{4 R^2 \varphi_1^2 \varphi_2^2}{(t_1+t_2)^2} ,} \qquad \displaystyle{\psi-8s\mu=\frac{4 R^2 \psi_1^2 \psi_2^2}{(t_1+t_2)^2} ,}
\end{array}
\end{equation}
where
\begin{equation}\notag
\begin{array}{l}
\varphi_1=\sqrt{2 s\sqrt{\tau}+ U_1}, \quad \varphi_2=\sqrt{2 s\sqrt{\tau}+ U_2}, \\
\psi_1=\sqrt{2 s\sqrt{\tau}- U_1}, \quad \psi_2=\sqrt{2 s\sqrt{\tau}- U_2}.
\end{array}
\end{equation}
This notation is intermediate, the signs chosen here do not affect the final expressions. The only important value is the new radical
\begin{equation}\label{neq4_32}
R=\sqrt{\vphantom{4s^2} t_1 t_2 + \sigma + U_1 U_2}.
\end{equation}
Now the expressions for $\mu_1,\mu_2$ become
\begin{equation}\label{neq4_33}
\begin{array}{l}
\displaystyle{\mu_1=\frac{R^2(\varphi_1 \varphi_2 - \psi_1 \psi_2)^2}{4 s (t_1+t_2)^2},} \qquad \displaystyle{\mu_2=\frac{R^2(\varphi_1 \varphi_2 + \psi_1 \psi_2)^2}{4 s (t_1+t_2)^2}. }
\end{array}
\end{equation}
We can avoid ``double roots'' by using the expanded form:
\begin{equation}\label{neq4_34}
\begin{array}{l}
\displaystyle{\mu_1=\frac{R^2(4 s^2 \tau + U_1 U_2- V_1 V_2)}{2 s (t_1+t_2)^2} ,} \qquad \displaystyle{\mu_2=\frac{R^2(4 s^2 \tau + U_1 U_2 + V_1 V_2)}{2 s (t_1+t_2)^2}.}
\end{array}
\end{equation}
Here we introduce new algebraic radicals
\begin{equation}
\label{neq4_35} V_1 = \sqrt{\mathstrut 4s^2\chi^2-t_1^2}, \quad V_2 = \sqrt{\mathstrut 4s^2\chi^2-t_2^2}
\end{equation}
having arbitrary signs.

To find $x_1,x_2$, in addition to (\ref{neq4_28}), (\ref{neq4_32}), (\ref{neq4_35}) denote
\begin{equation}
\label{neq4_36}
\begin{array}{ll} M_1=\sqrt{\vphantom{r^2} t_1+\tau+r^2}, & M_2=\sqrt{\mstrut t_2+\tau+r^2}, \\[1.5mm] N_1=\sqrt{\vphantom{r^2} t_1+\tau-r^2}, & N_2=\sqrt{\mstrut t_2+\tau-r^2}. \end{array}
\end{equation}
The signs here are also arbitrary. For the polynomials (\ref{neq4_17}) we have
\begin{equation}\label{neq4_37}
\begin{array}{l}
\displaystyle{\Phi_1=\frac{2 R^2 M_1^2 M_2^2}{(t_1+t_2)^2},} \qquad \displaystyle{\Phi_2=\frac{2 R^2 N_1^2 N_2^2}{(t_1+t_2)^2}.}
\end{array}
\end{equation}
Let $X=4r^4(x^2-\tau)+\tau(\Phi_1+\Phi_2)$. Then from (\ref{neq4_37}) and (\ref{neq4_30}) we find
\begin{equation}\notag
X^2-4 \tau^2 \Phi_1\Phi_2 = \frac{16 \tau^2 r^4(t_1-t_2)^2 R^4}{(t_1+t_2)^4}.
\end{equation}
Hence $2 \tau \sqrt{\mathstrut \Phi_1\Phi_2} = \sqrt{\mathstrut X_1 X_2}$, where
\begin{equation}\notag
\begin{array}{l}
\displaystyle{ X_1 = X+\sqrt{X^2-4 \tau^2 \Phi_1\Phi_2}=\frac{4\tau R^2 N_1^2 M_2^2}{(t_1+t_2)^2} ,}\\
\displaystyle{ X_2 = X-\sqrt{X^2-4 \tau^2 \Phi_1\Phi_2}=\frac{4\tau R^2 M_2^2 N_1^2}{(t_1+t_2)^2}.}
\end{array}
\end{equation}
For the numerators in (\ref{neq4_15}) we obtain
\begin{equation}\notag
4r^4(x^2-\tau) + \tau (\sqrt{\Phi_1} \pm \sqrt{\Phi_2})^2 = \frac{1}{2}(\sqrt{X_1} \pm \sqrt {X_2})^2.
\end{equation}
Then from (\ref{neq4_14}), (\ref{neq4_15}), (\ref{neq4_33}) we get
\begin{equation}\label{neq4_38}
\begin{array}{l}
\displaystyle{ x_1 = \frac{2 s \tau}{r^2}\left( \frac{M_2 N_1 + M_1 N_2} {\varphi_1 \varphi_2 + \psi_1 \psi_2} \right)^2,} \qquad \displaystyle{ x_2 = \frac{2 s \tau}{r^2}\left( \frac{M_2 N_1 - M_1 N_2} {\varphi_1 \varphi_2 - \psi_1 \psi_2} \right)^2,}
\end{array}
\end{equation}
or in the expanded form
\begin{equation}\label{neq4_39}
\begin{array}{l}
\displaystyle{x_1 = \frac{2 s \tau}{r^2} \frac{ (t_1+\tau)(t_2+\tau) - r^4 + M_1 N_1 M_2 N_2 }{ 4 s^2 \tau + U_1 U_2 + V_1 V_2},} \\[3mm]
\displaystyle{x_2 = \frac{2 s \tau}{r^2} \frac{ (t_1+\tau)(t_2+\tau) - r^4 - M_1 N_1 M_2 N_2 }{ 4 s^2 \tau + U_1 U_2 - V_1 V_2}.} \end{array}
\end{equation}

The expressions for $y_1,y_2$ in terms of $t_1,t_2$ are now easily found from the definition (\ref{neq4_6}) of $\mu_1, \mu_2$. Write
\begin{equation}\notag
\tau y_1 \mu_2 = r^2 x_1 \mu_2 - \mu^2, \quad \tau y_2 \mu_1 = r^2 x_2 \mu_1 - \mu^2.
\end{equation}
Then after substitution of (\ref{neq4_34}), (\ref{neq4_39}), (\ref{neq4_27}) we get
\begin{equation}\label{neq4_40}
\begin{array}{l}
\displaystyle{y_1  = 2 s \frac{\tau ( t_1 + t_2 - 2 p^2 + 2 \tau) - U_1 U_2 + M_1 N_1 M_2 N_2} {4 s^2 \tau + U_1 U_2 + V_1 V_2 },}\\[3mm]
\displaystyle{y_2  = 2 s \frac{\tau ( t_1 + t_2 - 2 p^2 + 2 \tau) - U_1 U_2 - M_1 N_1 M_2 N_2} {4 s^2 \tau + U_1 U_2 - V_1 V_2 }.}
\end{array}
\end{equation}

The variables $z_1,z_2$ are expressed from (\ref{neq4_17}), (\ref{neq4_18}):
\begin{equation}\label{neq4_41}
\begin{array}{l}
\displaystyle{z_1  = \frac{R} {\sqrt{2} \, r} \frac{M_1 M_2 + N_1 N_2}{t_1 + t_2} ,} \qquad \displaystyle{z_2  = \frac{R} {\sqrt{2} \, r} \frac{M_1 M_2 - N_1 N_2}{t_1 + t_2}.}
\end{array}
\end{equation}

To find the components of the angular velocity, start with $w_3$. Using the identity (\ref{neq4_30}) represent the polynomials $P,Q$ as follows
\begin{equation}
\notag \displaystyle{P=\frac{4\xi^2}{(t_1+t_2)^2}\,\widetilde{P}, \quad Q=\frac{2 \xi}{t_1+t_2}\,\widetilde{Q}.}
\end{equation}
Here
\begin{equation}\notag
\begin{array}{l}
\widetilde{P}=-(t_1^2-r^4)(t_2^2-r^4)+\tau[(2s^2-t_2)t_1^2+(2s^2-t_1)t_2^2-p^2(t_1^2+t_2^2)\\
\qquad \qquad + r^4(t_1+t_2-4s^2+2p^2)]+2\tau^2(2s^2-p^2)(t_1+t_2)\\
\qquad \qquad + \tau^3(t_1+t_2+4s^2-2p^2)+\tau^4,\\
\widetilde{Q}=t_1 t_2 +(2s^2-p^2)(t_1+t_2)+r^4+\tau(t_1+t_2+4s^2-2p^2)+\tau^2.
\end{array}
\end{equation}
In virtue of (\ref{neq4_35}), (\ref{neq4_36}) we obtain for the functions (\ref{neq4_25})
\begin{equation}
\notag
\begin{array}{l} \displaystyle{P_1  =  \frac{4\xi^2}{(t_1+t_2)^2} M_1^2 N_1^2 V_2^2,} \qquad  \displaystyle{P_2 =\frac{4\xi^2}{(t_1+t_2)^2} M_2^2 N_2^2 V_1^2.} \end{array}
\end{equation}
Then from (\ref{neq4_24}) using the identity (\ref{neq4_29}) we find
\begin{equation}
\label{neq4_42} \displaystyle{w_3=\frac{U_1-U_2}{\sqrt{2s\tau}}\frac{M_2 N_2 V_1 - M_1 N_1 V_2}{t_1^2-t_2^2}= \frac{1}{\sqrt{2s\tau}}\frac{M_2 N_2 V_1 - M_1 N_1 V_2}{U_1+U_2}}.
\end{equation}
Note that the expressions for $x_i,y_i,z_i$ contain only the product $V_1 V_2$. Therefore we can write the signs of the terms in the numerator of $w_3$ in different ways. Then, having fixed the form (\ref{neq4_42}) and applying (\ref{neq4_21}) to find $w_1,w_2$, we must point out the rule to define the signs of the radicals $\sqrt{\Theta_1},\sqrt{\Theta_2}$ in a way that guarantees all combinations satisfying the system of three equations (\ref{neq4_2}) linear in $w_i$, together with (\ref{neq4_42}). Since the determinant of this system is zero by virtue of (\ref{neq4_3}), we need to check only two of these equations, e.g.,
\begin{equation}\label{neq4_43}
\begin{array}{l}
(y_2 + 2s) w_1  + x_1 w_2  + z_1 w_3  = 0,  \\
x_2 w_1  + (y_1  + 2s) w_2  + z_2 w_3  = 0.
\end{array}
\end{equation}
From (\ref{neq4_38}) we write
\begin{equation}\notag
\begin{array}{l}
\displaystyle{ \sqrt{x_1} = \frac{\sqrt{2 s \tau}}{r} \frac{M_2 N_1 + M_1 N_2} {\varphi_1 \varphi_2 + \psi_1 \psi_2},} \qquad \displaystyle{ \sqrt{x_2} = \frac{\sqrt{2 s \tau}}{r} \frac{M_2 N_1 - M_1 N_2} {\varphi_1 \varphi_2 - \psi_1 \psi_2}.} \end{array}
\end{equation}
Here the signs are chosen to satisfy the necessary condition $\sqrt{x_1} \sqrt{x_2} \equiv x$, where the value $x$ is defined according to (\ref{neq4_27}). For $\Theta_1, \Theta_2$ the equations (\ref{neq4_22}) and (\ref{neq4_27}) then give
\begin{equation}\notag
\begin{array}{c}
(t_1+t_2)^2 \Theta_1  = -2 R^2 \varphi_1^2 \psi_2^2, \qquad (t_1+t_2)^2 \Theta_2  = -2 R^2 \psi_1^2 \varphi_2^2.
\end{array}
\end{equation}
Hence, introducing $\varepsilon_{1,2}=\pm 1$ we write
\begin{equation}\notag
\begin{array}{c}
(t_1+t_2) \sqrt{\mathstrut \Theta_1} = \ri \, \varepsilon_1 \, \sqrt{2} R \varphi_1 \psi_2, \qquad (t_1+t_2) \sqrt{\mathstrut \Theta_2} = \ri \, \varepsilon_2 \, \sqrt{2} R  \psi_1 \varphi_2.
\end{array}
\end{equation}
Then
\begin{equation}
\notag
\begin{array}{l} w_1=\displaystyle{\frac{r R}{4 s \sqrt{s \tau}} \frac {(\varepsilon_2 \varphi_2^2-\varepsilon_1 \psi_2^2)V_1+ (\varepsilon_1 \varphi_1^2-\varepsilon_2 \psi_1^2)V_2} {(t_1+t_2)(M_1 N_2 - M_2 N_1)}}, \\[3mm] w_2=\displaystyle{\frac{r R}{4 s \sqrt{s \tau}} \frac {(\varepsilon_2 \varphi_2^2-\varepsilon_1 \psi_2^2)V_1 - (\varepsilon_1 \varphi_1^2-\varepsilon_2 \psi_1^2)V_2} {(t_1+t_2)(M_1 N_2 + M_2 N_1)}}. \end{array}
\end{equation}
Choose $\varepsilon_2 = -\varepsilon_1$. The expressions
\begin{equation}
\notag w_1 = \pm \frac{r R}{\sqrt{s}(t_1+t_2)} \frac{V_1-V_2}{ M_2 N_1 - M_1 N_2}, \qquad w_2 = \mp \frac{r R}{\sqrt{s}(t_1+t_2)} \frac{V_1+V_2}{M_2 N_1 + M_1 N_2},
\end{equation}
do not satisfy (\ref{neq4_43}) for any choice of the signs. For $\varepsilon_1 = \varepsilon_2=\varepsilon$ out of two possibilities $\varepsilon=\pm 1$ only one satisfies (\ref{neq4_42}) and (\ref{neq4_43}), namely,
\begin{equation}
\label{neq4_44}
\begin{array}{l} \displaystyle{w_1=\frac{r (U_1 V_2+U_2 V_1  ) R} {2s \sqrt{s \tau}(t_1+t_2)(M_2 N_1 - M_1 N_2 )}},\\ \displaystyle{w_2=\frac{r (U_1 V_2-U_2 V_1) R} {2s \sqrt{s \tau}(t_1+t_2)(M_2 N_1 + M_1 N_2)}}. \end{array}
\end{equation}

To accomplish the algebraic solution note that we still have a ``double root'' expression for $R$. In fact, it can be eliminated. Put
\begin{equation}
\notag \varkappa = \sqrt{\sigma}.
\end{equation}
This value can be either real or pure imaginary. Denote
\begin{equation}
\notag
\begin{array}{l} K_1=\sqrt{\mstrut t_1+\varkappa}, \quad  K_2=\sqrt{\mstrut t_2+\varkappa},\quad  L_1=\sqrt{\mstrut t_1-\varkappa}, \quad L_2=\sqrt{\mstrut t_2-\varkappa}. \end{array}
\end{equation}
Then
\begin{equation}
\notag
\begin{array}{l} R=\ds{\frac{1}{\sqrt{\mathstrut 2}}}(K_1 K_2 + L_1 L_2),\qquad  U_1= K_1 L_1, \qquad U_2=K_2 L_2. \end{array}
\end{equation}

Finally, the formulae (\ref{neq4_39}), (\ref{neq4_40}), (\ref{neq4_41}), (\ref{neq4_42}) and (\ref{neq4_44}) algebraically define the values of complex phase coordinates (\ref{neq2_19}) for given $t_1,t_2$ up to the choice of signs for the following radicals
\begin{equation}\notag
\sqrt{s \tau},\, K_1,\,K_2,\, L_1,\, L_2, \, V_1,\, V_2,\, M_1,\, M_2,\, N_1,\, N_2.
\end{equation}
Here $V_i, M_i, N_i$ can be real or pure imaginary, and the pairs $K_1,K_2$ and $L_1,L_2$ can be real or complex conjugate. Accessible regions for $t_1,t_2$ are defined by the condition that the initial phase variables $\alpha_j,\beta_j,\omega_j$ $(j=1,2,3)$ in (\ref{neq2_19}) are real.

\subsection{The equations of motion and the real solution} To derive the differential equations, which describe the time evolution of the variables $t_1,t_2$ let us use the variables $s_1,s_2$ as the intermediate ones. From (\ref{neq2_41}), (\ref{neq4_7}) we get
\begin{equation}\label{neq4_48}
\displaystyle{s_1=\frac{\xi+\tau+r^2}{2x},}\quad \displaystyle{s_2=\frac{\xi+\tau-r^2}{2x}.}
\end{equation}
Calculate the derivatives by virtue of (\ref{neq3_3}). Then from (\ref{neq3_4}) we obtain
\begin{equation}\label{neq4_49}
\begin{array}{l}
\displaystyle{\frac{d s_1}{dt} = i \frac{r^2}{4 x^3}(z_1+z_2)(x_1 w_2-x_2 w_1),}\\[3mm]
\displaystyle{\frac{d s_2}{dt} = i \frac{r^2}{4 x^3}(z_1-z_2)(x_1 w_2+x_2 w_1).}
\end{array}
\end{equation}
On the other hand,
\begin{equation}
\label{neq4_50}
\begin{array}{l} \displaystyle{\frac{d s_1}{dt}=\frac{\partial s_1}{\partial t_1}\frac{d t_1}{dt}+\frac{\partial s_1}{\partial t_2}\frac{d t_2}{dt},}\quad \displaystyle{\frac{d s_2}{dt}=\frac{\partial s_2}{\partial t_1}\frac{d t_1}{dt}+\frac{\partial s_2}{\partial t_2}\frac{d t_2}{dt}.} \end{array}
\end{equation}
Substitute (\ref{neq4_27}) to (\ref{neq4_48}) to obtain
\begin{equation}
\label{neq4_51}
\begin{array}{ll} \displaystyle{\frac{\partial s_1}{\partial t_1} = - \frac{t_1 t_2-\sigma +U_1 U_2}{2\sqrt{\tau}(t_1-t_2)^2} \frac {M_2^2}{U_1},} & \displaystyle{\frac{\partial s_1}{\partial t_2} = \frac{t_1 t_2-\sigma+ U_1 U_2}{2\sqrt{\tau}(t_1-t_2)^2} \frac {M_1^2}{U_2},} \\[4mm] \displaystyle{\frac{\partial s_2}{\partial t_1} = - \frac{t_1 t_2-\sigma+U_1 U_2}{2\sqrt{\tau}(t_1-t_2)^2} \frac {N_2^2}{U_1},} & \displaystyle{\frac{\partial s_2}{\partial t_2} = \frac{t_1 t_2-\sigma+ U_1 U_2}{2\sqrt{\tau}(t_1-t_2)^2} \frac {N_1^2}{U_2}.} \end{array}
\end{equation}
From (\ref{neq4_50}), (\ref{neq4_51}) we have
\begin{equation}
\label{neq4_52}
\begin{array}{ll} \displaystyle{\frac{d t_1}{d t} = \frac{\sqrt{\tau}(t_1-t_2) U_1}{r^2(t_1 t_2-\sigma + U_1 U_2)} ( M_1^2 \frac{ds_2}{dt} -  N_1^2 \frac{ds_1}{dt}),} \\[4mm] \displaystyle{\frac{d t_2}{d t} = \frac{\sqrt{\tau}(t_1-t_2) U_2}{r^2( t_1 t_2-\sigma + U_1 U_2)} (M_2^2 \frac{ds_2}{dt} - N_2^2 \frac{ds_1}{dt} ).} \end{array}
\end{equation}
Express the values (\ref{neq4_49}) in terms of $t_1,t_2$ using (\ref{neq4_39}), (\ref{neq4_41}), (\ref{neq4_44}):
\begin{equation}
\notag
\begin{array}{ll} \displaystyle{\frac{d s_1}{d t} = \ri \frac{(t_1 t_2+\sigma+U_1 U_2) M_1 M_2(M_1 N_2 V_2 - M_2 N_1 V_1)}{2\sqrt{2 s}\,\tau (U_1-U_2)^2(t_1-t_2)},}\\[3mm] \displaystyle{\frac{d s_2}{d t} = \ri \frac{(t_1 t_2+\sigma+U_1 U_2) N_1 N_2( M_2 N_1 V_2 - M_1 N_2 V_1)}{2\sqrt{2 s}\,\tau (U_1-U_2)^2(t_1-t_2)}.} \end{array}
\end{equation}
Substitute these expressions to (\ref{neq4_52}). Then in view of the relation (\ref{neq4_31}) after some obvious transformations we come to the equations of the Kowalevski type
\begin{equation}
\notag
\begin{array}{ll} \displaystyle{(t_1-t_2)\frac{d t_1}{d t} = \sqrt{\frac{1}{2 s \tau}(4s^2\chi^2-t_1^2)(t_1^2-\sigma)[r^4-(t_1+\tau)^2]}\, ,} \\ [3mm] \displaystyle{(t_1-t_2)\frac{d t_2}{d t} = \sqrt{\frac{1}{2 s \tau}(4s^2\chi^2-t_2^2)(t_2^2-\sigma)[r^4-(t_2+\tau)^2]}\,.} \end{array}
\end{equation}

To obtain the expressions for the real phase variables use the inverse linear mapping to the change~(\ref{neq2_19}). From (\ref{neq4_39}), (\ref{neq4_40}), (\ref{neq4_41}) we find $\alpha_j,\beta_j$ $(j=1,2,3)$ as functions of $t_1,t_2$:
\begin{equation}
\begin{array}{l} \displaystyle{\alpha_1=\frac{(\mathcal{A}-r^2 U_1 U_2)(4 s^2 \tau+U_1 U_2)-(\tau+r^2) M_1 N_1 M_2 N_2 V_1 V_2}{4 r^2 s\, \tau (U_1+U_2)^2}, } \\[3mm] \displaystyle{\alpha_2=\ri \frac{(\mathcal{A} -r^2 U_1 U_2)V_1 V_2 -(4 s^2 \tau+U_1 U_2)(\tau+r^2) M_1 N_1 M_2 N_2}{4 r^2 s\, \tau (U_1+U_2)^2}, }\\[3mm] \displaystyle{\alpha_3= \frac{R }{r \sqrt{\mathstrut 2} } \, \frac {M_1 M_2}{t_1+t_2},} \\ \displaystyle{\beta_1=\ri \frac{(\mathcal{B}+r^2 U_1 U_2)V_1 V_2-(4 s^2 \tau+U_1 U_2)(\tau-r^2) M_1 N_1 M_2 N_2}{4 r^2 s\, \tau (U_1+U_2)^2}, } \\[3mm] \displaystyle{\beta_2=-\frac{(\mathcal{B} + r^2 U_1 U_2)(4 s^2 \tau + U_1 U_2)-(\tau-r^2) M_1 N_1 M_2 N_2 V_1 V_2}{4 r^2 s\, \tau (U_1+U_2)^2}, }\\[3mm] \displaystyle{\beta_3= - \ri \frac{R}{r \sqrt{\mathstrut 2} } \, \frac { N_1 N_2}{t_1+t_2}.} \end{array}\label{neq4_54}
\end{equation}
Here for brevity we put
\begin{equation}\notag
\begin{array}{l} \mathcal{A}=[(t_1+\tau+r^2)(t_2+\tau+r^2)-2(p^2+r^2)r^2]\tau, \\ \mathcal{B}=[(t_1+\tau-r^2)(t_2+\tau-r^2)+2(p^2-r^2)r^2]\tau. \end{array}
\end{equation}

The angular velocities $\omega_j$ $(j=1,2,3)$ are found from (\ref{neq2_19}), (\ref{neq4_42}), (\ref{neq4_44}):
\begin{equation}
\label{neq4_56}
\begin{array}{l} \displaystyle{\omega_1=  \frac{R }{4 r s\, \sqrt{s \,\tau}} \, \frac{M_2 N_1 U_1 V_2 + M_1 N_2 U_2 V_1}{ t_1^2-t_2^2},} \\[3mm] \displaystyle{\omega_2= -\frac{\ri\,R}{4 r s\, \sqrt{s \,\tau} } \, \frac{M_2 N_1 U_2 V_1 + M_1 N_2 U_1 V_2}{t_1^2-t_2^2},} \\[3mm] \displaystyle{\omega_3=  \frac{U_1-U_2}{\sqrt{2s\tau}}\frac{M_2 N_2 V_1 - M_1 N_1 V_2}{t_1^2-t_2^2}.} \end{array}
\end{equation}
 These explicit expressions of the phase variables in terms of the $t_1,t_2$ along with the separated differential equations give the complete analytical solution for the third critical subsystem.

\section{Conclusion}
In this article we presented the results dealing with analytical solutions in the generalized Kowalevski problem. The questions of qualitative and topological analysis of the arising subsystems are connected with the investigation of the multi-valued dependencies (\ref{neq3_15}) and (\ref{neq4_54}), (\ref{neq4_56}). Formally these dependencies are the inverse mappings for the coverings of degree $2^n$, where $n$ is the number of radicals with arbitrary signs. Lately, new methods of investigating the phase topology of algebraically solvable systems were proposed in \cite{KhRCD09} based on calculating the invariants of some $\mathbb{Z}_2$-linear mappings. (When this article was already published a theory of the topological analysis based on Boolean vector-functions was presented in \cite{KhND10,KhND11} with applications to the classical problems and to the systems described above.)

Up to this moment the first critical subsystem (the Bogoyavlensky case) has not received any analytical solution despite the fact that its topology is completely investigated in \cite{ZotRCD} and various Lax representations are found in \cite{BogEn,ReySem}. For the case of the Kowalevski gyrostat in two constant fields no analytical solutions are found except for the special periodic motions and their bifurcations studied in \cite{Kh37}, though all critical subsystems are analytically described in \cite{KhND07}.

The Kowalevski problem in the rigid body dynamics remains as challenging as it was during the last 120 years.

The Russian version of this article was written with the financial support of the RFBR and the Volgograd Region Administration (grant No\,10-01-97001).

\end{document}